\documentclass[11pt,a4paper]{amsart}
\usepackage[T1]{fontenc}
\usepackage{lmodern}
\usepackage[a4paper,margin=28mm,headheight=15pt]{geometry}
\usepackage{amsmath,amssymb,amsthm,mathtools}
\usepackage{mathrsfs}
\usepackage{microtype}
\usepackage{enumitem}
\usepackage{booktabs,array}
\usepackage{xcolor}

\definecolor{linkblue}{RGB}{0,0,0}
\usepackage[colorlinks=true,linkcolor=linkblue,citecolor=linkblue,urlcolor=linkblue]{hyperref}
\usepackage[nameinlink,noabbrev]{cleveref}
\usepackage{fancyhdr}
\allowdisplaybreaks[2]
\setlist[enumerate]{itemsep=0.25em,topsep=0.5em}
\setlist[itemize]{itemsep=0.2em,topsep=0.4em}
\numberwithin{equation}{section}
\newtheorem{theorem}{Theorem}[section]
\newtheorem{proposition}[theorem]{Proposition}
\newtheorem{lemma}[theorem]{Lemma}
\newtheorem{corollary}[theorem]{Corollary}
\theoremstyle{definition}
\newtheorem{definition}[theorem]{Definition}

\theoremstyle{remark}
\newtheorem{remark}[theorem]{Remark}
\newcommand{\ddc}{dd^c}
\newcommand{\dbar}{\bar\partial}
\newcommand{\I}{\mathcal I}
\newcommand{\OO}{\mathcal O}
\newcommand{\HH}{\mathcal H}

\newcommand{\C}{\mathbb C}
\newcommand{\R}{\mathbb R}
\newcommand{\Z}{\mathbb Z}
\newcommand{\one}{\mathbf 1}
\DeclareMathOperator{\Dom}{Dom}
\DeclareMathOperator{\Ran}{Ran}
\DeclareMathOperator{\Spec}{Spec}
\DeclareMathOperator{\rank}{rank}

\DeclareMathOperator{\PSH}{PSH}
\DeclareMathOperator{\supp}{supp}
\DeclareMathOperator{\graph}{graph}

\title{Cohomology growth for pseudoeffective line bundles on compact K\"ahler manifolds}
\author[Hongrong Chen]{Hongrong Chen}
\address{Hongrong Chen: School of Mathematics (Zhuhai), Sun Yat-sen University, Zhuhai 519082, China.}
\email{chenhr73@mail2.sysu.edu.cn}

\author[Guokuan Shao]{Guokuan Shao}
\address{Guokuan Shao: School of Mathematics (Zhuhai), Sun Yat-sen University, Zhuhai 519082, China.}
\email{shaogk@mail.sysu.edu.cn}
\thanks{The second author was supported by National Key R\&D Program of China No. 2024YFA1015200 and by National Natural Science Foundation of China No. 12471082.}

\author[Wenxuan Wang]{Wenxuan Wang}
\address{Wenxuan Wang: Institute of Mathematics, Academy of Mathematics and Systems Science, Chinese Academy of Sciences, Beijing 100190, China.}
\email{wxwang@amss.ac.cn}

\author[Xiangyu Zhou]{Xiangyu Zhou}
\address{Xiangyu Zhou: Institute of Mathematics, Academy of Mathematics and Systems Science, Chinese Academy of Sciences, Beijing 100190, China.}
\email{xyzhou@math.ac.cn}
\thanks{The fourth author was supported by National Key R\&D Program of China No. 2021YFA1003100 and by National Natural Science Foundation of China No. 12288201}

\subjclass[2020]{Primary 32L10; Secondary 32C35}
\date{}
\begin{document}

\begin{abstract}
Let $L$ be a pseudoeffective line bundle with a singular metric $h$ satisfying $\Theta_h(L)\geq0$ on a compact K\"ahler manifold $X$ of dimension $n$. Let $\I(h^k)$ denote the multiplier ideal sheaf of $h^k$. For every holomorphic vector bundle $E$ and $0\le q\le n$, we prove that
\[
 \dim_{\C}H^q\!\left(X,E\otimes L^k\otimes\I(h^k)\right)
 =O(k^{n-q}),
\]
which gives an affirmative answer to Matsumura's growth problem. As an application, we resolve the cohomology growth problem for nef line bundles raised by Demailly.
\end{abstract}
\maketitle
\thispagestyle{plain}
\medskip
\noindent\textbf{Keywords.} Singular Hermitian metric, multiplier ideal sheaf, pseudoeffective line bundle, nef line bundle, Dolbeault cohomology, spectral localization.

\section{Introduction}
\subsection{Solutions to Demailly’s and Matsumura’s Problems}
Let $L$ be a holomorphic line bundle on a compact K\"ahler manifold
$X$ of dimension $n$, and let $E$ be a fixed holomorphic vector
bundle. If $L$ is positive, $H^q(X, E\otimes L^k)=0$
for every $q \geq 1$ and all sufficiently large $k$. When the curvature of $L$ is weakly positive, one may then ask how rapidly its dimension can grow. For a semipositive line bundle, Berndtsson \cite{Berndtsson} proved $\operatorname{dim}_{\C} H^q(X,E\otimes L^k)=O(k^{n-q})$ for any compact complex manifolds.
For a nef line bundle on a projective manifold, Demailly established a similar result and asked whether the estimate extends beyond the algebraic case \cite[Lemma~6.18]{Demailly}. Matsumura proved the projective analog for the pseudoeffective line bundle with singular semipositive metrics $h$, with the multiplier ideal $\I(h^k)$ in the coefficient sheaf \cite[Theorem~4.1]{MatsumuraNadel}, and formulated a compact K\"ahler extension problem in \cite[Problem~3.4]{MatsumuraSurvey}. In this paper, we resolve the two problems.

\begin{theorem}\label{thm1}
Let $(X,\omega)$ be a compact K\"ahler manifold of complex dimension $n$, let $(L,h)$ be a holomorphic line bundle with a singular Hermitian metric whose Chern curvature current $\Theta_h(L)$ satisfies $\Theta_h(L)\ge0$, and let $E$ be a fixed holomorphic vector bundle. Then
\begin{equation}
 \operatorname{dim}_{\C}H^q\!\left(X,E\otimes L^k\otimes\I(h^k)\right)
 =O(k^{n-q}),\quad 0\le q\le n.
\end{equation}
\end{theorem}

The case $q=0$ follows from the general estimate for holomorphic sections \cite[Theorem~5.1]{DemaillyMorse}.  As an application, by using a uniform version of the argument for smooth metrics, we have the following

\begin{corollary}\label{cor1}
Let $X$ be a compact K\"ahler manifold of dimension $n$, $L$ a nef holomorphic line bundle, and $E$ a fixed holomorphic vector bundle. Then
\begin{equation}
\operatorname{dim}_{\C}H^q(X,E\otimes L^k)=O(k^{n-q}),\quad 0\le q\le n.
\end{equation}
\end{corollary}

These results answer the compact K\"ahler cases of the questions raised by Matsumura and Demailly. After obtaining our results, we became aware of the recent work of Qin \cite[Theorem~B]{Qin}, which obtains also a similar result of Corollary \ref{cor1}.

\subsection{Background}
The asymptotic estimates for smooth semipositive line bundles were obtained by Siu \cite{SiuVanishing} and by Demailly's holomorphic Morse inequalities \cite{DemaillyMorse}. For smooth semipositive metrics and line-bundle coefficients, Berndtsson proved a stronger low-eigenvalue estimate \cite[Theorem~1.1]{Berndtsson}. Z. Wang and the fourth named author  treated arbitrary vector-bundle coefficients in the smooth case and algebraic singularities in the singular case \cite[Theorems~1.3-1.4 and Questions~1.1-1.2]{WZ}. Their treatment of algebraic singularities combines the smooth estimate with the modification-based approach of Bonavero's singular holomorphic Morse inequalities \cite{BonaveroMorse}. For line-bundle coefficients, see also the work of H.~Wang \cite{HuanWangGrowth}. Our main theorem treats arbitrary plurisubharmonic singularities on compact K\"ahler manifolds.

A general framework for this question is provided by Demailly's holomorphic Morse inequalities \cite{DemaillyMorse}. They relate the asymptotic dimensions of Dolbeault cohomology to the signature of the curvature, which yield $\dim_{\C}H^q(X,E\otimes L^k)=o(k^n)$ for $q \geq 1$, when $L$ is semipositive. This provides motivation for seeking more refined estimates, for example, estimates of order $O(k^{n-q})$ considered in this paper.

Holomorphic Morse inequalities have been developed in several directions. Bismut \cite{BismutMorse} gave a heat equation proof by expressing the relevant inequalities in terms of alternating heat traces. In the singular setting, Bonavero \cite{BonaveroMorse} established holomorphic Morse inequalities involving multiplier ideals for metrics with analytic singularities. In non-compact setting,  Bouche \cite{BoucheNoncompactMorse} treated noncompact $q$-convex manifolds under suitable hypotheses, while Berman \cite{BermanLocalMorse} obtained local holomorphic Morse inequalities from pointwise Bergman kernel estimates, with no compactness assumption. More recently, M.~Liu, together with the second and third authors \cite{LiuShaoWenxuanWangStrong}, provided a unified treatment of strong holomorphic Morse inequalities in several geometric settings under suitable analytic and curvature assumptions. For CR manifolds, C.-Y. Hsiao, R.-T. Huang, X. Li, and the second author \cite{HsiaoHuangLiShaoBoundary} established $S^1$-equivariant index formulas and Morse inequalities under a compatible holomorphic $S^1$-action transversal to the boundary. In \cite[Theorem~1.6]{HsiaoHuangLiShaoBoundary}, they also obtained sharp growth estimates for equivariant holomorphic functions when the boundary is weakly pseudoconvex everywhere and strongly pseudoconvex at one point. Through the disc-bundle correspondence, this result relates to the cohomology growth problem when $q=0$ in the smooth metric setting. Z. Wang and the fourth author \cite{WangZhouJDG} established sharp spectral estimates for Fourier components of the Kohn Laplacian on compact weakly pseudoconvex CR manifolds with transversal $S^1$-actions, yielding degree-dependent growth bounds for Kohn--Rossi cohomology with applications to CR Morse inequalities.

\subsection{Outline of the proof}
We establish the same local mass estimate of order $r^{2q}$ for weighted harmonic forms in Section \ref{sec3} that appears in Berndtsson \cite{{Berndtsson}} and Z. Wang and the fourth named author \cite{WZ}. A further step is needed to deduce a dimension estimate. In the smooth setting, based on the local mass bound, Berndtsson obtained a pointwise estimate by rescaling the Laplacian on balls of radius $k^{-1/2}$. Z. Wang and the fourth named author extended the argument for bundle-valued version. Such pointwise estimate controls the Bergman kernel function and hence the dimension.

For an arbitrary plurisubharmonic weight, the curvature lower bound does not provide the local quadratic expansion of the weight and the uniform control of the rescaled operator (cf. the proof of Theorem 2.3 in \cite{{Berndtsson}}). Replacing the singular weight by a smooth approximation yields a weighted $\bar\partial$-Laplacian with smooth coefficients. However, the curvature lower bound alone does not ensure uniform bounds for the coefficients of the rescaled Laplacian. For metrics with algebraic singularities, Z. Wang and the fourth named author developed techniques from Bonavero \cite{BonaveroMorse}. They resolved the singularities and twisted a suitable power of the pulled-back line bundle by a divisor line bundle, which cancels the logarithmic part of the weight, producing a line bundle with a smooth semipositive metric to which their smooth estimates can apply \cite[Section~6]{WZ}. 

In our setting there is no guarantee that its pullback under a modification can be written locally as the sum of a smooth function and logarithmic terms along a divisor. We therefore develop a dimension count argument. We construct local Gaussian sections and a spectral projection $\Pi_{k,\ell}$  of rank $O(k^n)$ with the crucial estimate \eqref{eq54} in Section \ref{sec4}, which is needed in Proposition~\ref{prop3} for estimating the localized sections $\chi_i\sigma_i\gamma$. Combined with the local mass estimate in Corollary \ref{cor2}, we obtain an upper bound $Cr^{2q}\operatorname{rank}\Pi_{k,\ell}$ in \eqref{eq68}. Taking $r=k^{-1/2}$ yields the required dimension bound of order $O(k^{n-q})$.

The paper is organized as follows. In Section \ref{sec2} we introduce basic notions and technical lemmas needed. In Section \ref{sec3}, we prove the local mass estimates for weighted harmonic forms associated with singular semipositive metrics. In Section \ref{sec4}, we construct local Gaussian sections and spectral spaces, and establish the spectral estimates needed. Section \ref{sec5} is devoted to proving the main theorem and the corollary by a spectral localization approach. The appendices contain some supplementary results.

\subsection{AI Disclosure}
After developing our idea on spectral localization approach for estimating the dimension of cohomology groups, we were suggested by ChatGPT to consult Laeng's work \cite{LaengThesis}.

\section{Preliminary}\label{sec2}
Throughout this paper, we assume that $(X,\omega)$ is a compact K\"ahler manifold of complex dimension $n$ with a fixed K\"ahler form $\omega$. Let $0\le q\le n$, set $p=n-q$, we use the conventions
\[
d^c=\frac{i}{2}(\bar\partial-\partial),\quad
 \ddc =i\partial\bar\partial,\quad
 \omega_q=\frac{\omega^q}{q!},\quad dV=\omega_n.
\]
 We write $\PSH(X,\omega)$ for the upper semicontinuous, locally integrable functions $v$ satisfying $\omega+dd^cv\ge0$. Let $\beta=\frac{i}{2}\sum_{j=1}^n dz_j\wedge d\bar z_j$ be the standard K\"ahler form on $\mathbb C^n$. Let $E, F$ be holomorphic vector bundles with fixed smooth Hermitian metrics $h_E$, $h_F$ respectively over $X$. Denote by $K_X=\Lambda^{n,0}T^*X$ the canonical line bundle of $X$.

We write $I$ for the identity operator. For a self-adjoint operator $T$ and a Borel set $J\subset\mathbb R$, $\mathbf{1}_J(T)$ denotes the corresponding spectral projection. The constants in the proofs of the results, such as $c, C, C_A$, may vary, but the (in)dependence remains unchanged.

\subsection{Singular Hermitian metrics and multiplier ideals}
Let $(L,h)$ be a holomorphic line bundle over $X$ endowed with a singular Hermitian metric $h$. Denote by $\Theta_h(L)$ the Chern curvature current. In a local holomorphic frame $e_L$ of $L$, we write $|e_L|_{h}^2=e^{-\phi}$, where $\phi$ is the local weight of $h$.

\begin{definition}
The multiplier ideal sheaf associated with $h^k$ is
\[
 \I(h^k)_x=\left\{f\in\OO_{X,x}: |f|^2e^{-k\phi}\text{ is locally integrable near }x\right\}.
\]
\end{definition}

\begin{definition}
A holomorphic line bundle $L$ is called pseudoeffective if it admits a singular Hermitian metric $h$ with $\Theta_h(L)\geq 0$.
\end{definition}
In the sequel sections, we assume that $(L,h)$ is pseudoeffective with $\Theta_h(L)\ge 0$. The following is an analytic characterization of a nef line bundle \cite{DPSNef} .

\begin{definition}
A holomorphic line bundle $L$ is called nef if for every $\varepsilon>0$, there exists a smooth Hermitian metric $g_\varepsilon$ with $\Theta_{g_\varepsilon}(L)\geq-\varepsilon\omega$.
\end{definition}

For a smooth Hermitian connection $D$ on a holomorphic vector bundle $E$, write $R^D=D^2$ for its curvature and $\Theta_D(E)=iR^D$. On a line bundle $L$, $[\Theta_D(L)/(2\pi)]=c_1(L)$. For the Chern connection of $(L,h)$, $\Theta_h(L)=\ddc\phi$. Fix a smooth metric $h_0$ on $L$ and write $h=h_0e^{-U}$. Then $\Theta_{h_0}(L)+\ddc U\ge0$. We assume that $U\le0$.

\subsection{A smooth approximation for q.p.s.h. function}

We use the B\l ocki-Ko\l odziej regularization theorem \cite[Theorem~1]{BlockiKolodziej} with a normalization adapted to the weighted comparison below.

\begin{lemma}\label{lem1}
Let $U\in\PSH(X,\omega)$ satisfying $U\le0$. Then there exist smooth functions $U_\nu\le0$ and nonnegative numbers $a_\nu\to0$ such that
\begin{equation}\label{eq3}
 \begin{gathered}
 U_\nu\ge U-a_\nu,\quad
 \ddc U_\nu\ge-\omega,\quad
 U_\nu\rightarrow U \ \text{a.e.}
 \end{gathered}
\end{equation}
in $L^1(X, \omega_n)$.
\end{lemma}

\begin{proof}
By \cite[Theorem~1]{BlockiKolodziej}, there are smooth $V_\nu\in\PSH(X,\omega)$ decreasing pointwisely to $U$. Set
\[
 s_\nu=\sup_XV_\nu,\quad
 a_\nu=\max\{0,s_\nu\},\quad U_\nu=V_\nu-a_\nu.
\]
The sequence $s_\nu$ decreases to some $s\ge\sup_XU$. Choose maximum points $x_\nu$ which achieve maximal values and a subsequence converging to $x$. For each $j$, $s_\nu\le V_j(x_\nu)$ for all sufficiently large $\nu$ and then $s\le V_j(x)$. Letting $j\to\infty$ yields $s\le U(x)\le\sup_XU$. Thus $s=\sup_XU\le0$ and $a_\nu$ converges decreasingly to $0$.

It follows that $U_\nu\le0$, $U_\nu\ge U-a_\nu$, $\ddc U_\nu\ge-\omega$, and $U_\nu\to U$ pointwisely. Since $U\in L^1(X)$ and
\[
 0\le-U_\nu\le-U+a_\nu\le-U+a_1
\]
almost everywhere, the dominated convergence theorem shows $L^1$ convergence, which proves \eqref{eq3}.
\end{proof}

\subsection{The weighted Dolbeault complex and harmonic representatives}
Define
\[
L^2_{n,q}(X,F\otimes L^k,h^k)=\bigg\{u\in
L^2\big(X,\Lambda^{n,q}T^*X\otimes F\otimes L^k\big):\int_X|u|_{\omega,h_F,h^k}^2\,dV<+\infty\bigg\}.
\]
Note that the pointwise norm $|\cdot|_{\omega,h_F,h^k}$ is induced by the Kähler form $\omega$, the smooth Hermitian metric $h_F$ on $F$, and the singular Hermitian metric $h^k$ on $L^k$. We write $|\cdot|_{h^k}=|\cdot|_{\omega,h_F,h^k}$ for short without confusion. Let $\|u\|^2_{h^k}=\int_X|u|^2_{h^k}dV$. Define
\[
\Dom\dbar=\{u\in L^2_{n,q}(X,F\otimes L^k,h^k):
\dbar u\in L^2_{n,q+1}(X,F\otimes L^k,h^k)\text{ in the sense of distributions}\}.
\]
The maximal distributional realization of $\dbar$ is closed and densely defined (cf. \cite[Section 1.6]{Takegoshi1997}). Denote by $\dbar^*_{h^k}$ the Hilbert adjoint of $\dbar$ with respect to the metric induced by $\omega$, $h_F$ and $h^k$. Define
\[
\HH_k^{n,q}(F)=\ker\dbar\cap\ker\dbar^*_{h^k},
\]
which is the space of weighted harmonic $(n,q)$-forms with values in $F\otimes L^k$. Note that $u$ is not necessarily smooth.

\begin{lemma}\label{lem2}
The maximal $\dbar$ operators have closed range and
\begin{equation}
 \HH_k^{n,q}(F)\simeq
 H^q\!\left(X,K_X\otimes F\otimes L^k\otimes\I(h^k)\right).
\end{equation}
In particular, $\HH_k^{n,q}(F)$ is finite dimensional.
\end{lemma}

\begin{proof}
By the local weighted $L^2$ solvability of $\bar\partial$, the complex of $(\Dom\dbar,\dbar)$ provides a fine resolution of
\[
K_X\otimes F\otimes L^k
\otimes\I(h^k).
\]
Thus
\[
\frac{\ker\dbar}
{\operatorname{Im}\dbar}
\simeq
H^q\bigl(X,K_X\otimes F\otimes L^k
\otimes\I(h^k)\bigr).
\]
Hence, by the Cartan--Serre finiteness theorem,
\[
\dim_{\mathbb C}
\frac{\ker\dbar}
{\operatorname{Im}\dbar}
<\infty.
\]
Consider the operator
\[
\dbar:
\Dom\dbar\cap L^2_{n,q-1}(X,F\otimes L^k,h^k)
\longrightarrow \ker\dbar\cap L^2_{n,q}(X,F\otimes L^k,h^k).
\]
Equip $\Dom \dbar$ with the graph norm. This is a bounded operator between Hilbert spaces, and its range has finite codimension, which implies that $\operatorname{Im}\dbar$ is closed; see \cite[Proposition 6.2.12] {BogachevSmolyanov2020}.
Finally, using the closeness of image of $\dbar$, we get $\HH_k^{n,q}(F)\simeq H^q\!\left(X,K_X\otimes F\otimes L^k\otimes\I(h^k)\right)$ by the standard $L^2$ Hodge decomposition.
\end{proof}

\subsection{The Lefschetz transform and the smooth Siu inequality}
 Let $*_\omega$ denote the Hodge star operator with respect to $\omega$. We have $\left\langle u\wedge\overline{*_\omega u} \right\rangle_{h^k}=|u|_{h^k}^2\,dV$ for an $(n,q)$-form $u$ with values in $F\otimes L^k$. Note that if $u=\alpha\otimes e_F\otimes e_L^k$ locally, then $\left\langle u\wedge\overline{*_\omega u} \right\rangle_{h^k}=|e_F|^2_{h_F}e^{-k\phi}\alpha\wedge\overline{*_\omega\alpha}$. Define
\[\gamma=*_\omega u,
\]
which is a $(p,0)$-form with values in $F\otimes L^k$. It is known that \cite[Section~3]{Berndtsson}
\begin{equation}
|\gamma|_{h^k}=|u|_{h^k}, \quad u=i^{p^2}\gamma\wedge\omega_q.
\end{equation}
Define
\[
T_u=i^{p^2}\left\langle
\gamma\wedge\overline{\gamma}\right\rangle_{h^k}.
\]
Then we have
\begin{equation}
T_u\wedge\omega_q=|u|_{h^k}^2\,dV.
\end{equation}
For $q\ge1$, $T_u\wedge\omega_{q-1}$ is a positive $(n-1,n-1)$-current.

\begin{lemma}\label{lem3}
Suppose that $h$ is smooth with $\Theta_h(L)\ge 0$. Let $u$ be a smooth $(n,q)$-form with values in $F\otimes L^k$. Then there exists a constant $C\ge0$ depending only on $\omega, h_F$, such that every nonnegative $f\in C_c^\infty(X)$ satisfies
\begin{equation}\label{eq7}
 \int f|\dbar\gamma|_{h^k}^2\,dV\le \bigl\langle\ddc T_u\wedge\omega_{q-1},f\bigr\rangle
+C\int f|u|_{h^k}^2\,dV+\int f|\dbar u|_{h^k}^2\,dV
+2\bigl|\langle\dbar^*_{h^k}u,\dbar^*_{h^k}(fu)\rangle_{h^k}\bigr|,
\end{equation}
where $\langle\dbar^*_{h^k}u, \dbar^*_{h^k}(fu)\rangle_{h^k}=\int_X\langle\dbar^*_{h^k}u,\dbar^*_{h^k}(fu)\rangle_{h^k}\,dV$.
\end{lemma}

\begin{proof}
It follows from \cite[Lecture 3, Proposition~4.1]{BerndtssonLectures} that
\begin{align}\label{eq8}
 \ddc T_u\wedge\omega_{q-1}
=\bigl(&-2\operatorname{Re}\langle\dbar\dbar^*_{h^k}u,u\rangle_{h^k}
+\langle\Theta_{F\otimes L^k}\wedge\Lambda u,u\rangle_{h^k}
\notag\\
 &+|\dbar\gamma|_{h^k}^2+|\dbar^*_{h^k}u|_{h^k}^2-|\dbar u|_{h^k}^2\bigr)dV.
\end{align}
Since $\Theta_{F\otimes L^k}$ is smooth, \eqref{eq8} implies immediately that for some constant $C>0$
\[
\ddc T_u\wedge\omega_{q-1}
\ge\Bigl(|\bar\partial\gamma|_{h^k}^2-C|u|_{h^k}^2-|\bar\partial u|_{h^k}^2
-2\operatorname{Re}\langle\bar\partial\dbar^*_{h^k}u,u\rangle_{h^k}\Bigr)dV.
\]
Then \eqref{eq7} follows.
\end{proof}

\section{Local estimates in singular metrics}\label{sec3}

To obtain cohomology growth for semipositive line bundle, Berndtsson \cite{Berndtsson} developed the local mass estimate by using Siu's $\partial\bar\partial$-formula for $T_u$ to control $T_u\wedge\omega_q=|u|_{h^k}^2dV$ on balls with order $r^{2q}$. Z. Wang and the fourth named author \cite{WZ} established the Siu's formula for singular metrics and local mass estimates for holomorphic vector bundle in smooth semipositive metrics. We work on the case for weakly harmonic forms and arbitrary plurisubharmonic local weights.

\subsection{The Siu type inequality for singular metrics}
The key point in Proposition~\ref{prop1} is to show the $L^2$ integrablity of $\dbar\gamma$ by using local smooth approximation to $u$, $\dbar u$ and $\dbar^*_{h^k}u$ simultaneously.

\begin{proposition}\label{prop1}
For every $u\in\HH_k^{n,q}(F)$, we have $\dbar\gamma\in L^2_{p,1}(X,F\otimes L^k,h^k)$. There exists a constant $C\ge0$ depending only on $\omega, h_F$, such that
\begin{equation}\label{eq9}
 \ddc T_u\wedge\omega_{q-1}
 \ge |\dbar\gamma|_{h^k}^2\,dV
-C|u|_{h^k}^2\,dV
\end{equation}
in the sense of currents. In particular,
\begin{equation}\label{eq10}
 \|\dbar\gamma\|_{h^k}^2\le C\|u\|_{h^k}^2.
\end{equation}
\end{proposition}

\begin{proof}
Fix a local holomorphic frame $e_L$ with $|e_L|_h^2=e^{-\phi}$. Choose a relatively compact domain $D_2$ on which $L$ is trivialized. The local weight $k\phi$ of $h^k$ is plurisubharmonic. There exist smooth plurisubharmonic functions $\phi_\nu$ which converge decreasingly to $\phi$ on a neighborhood of $D_1\Subset D_2$. Then $0<c\le{e^{-k\phi_\nu}}\le{e^{-k\phi}}$ for some constant $c$ and ${e^{-k\phi_\nu}}$ converge increasing to ${e^{-k\phi}}$ almost everywhere on $D_1$. Here we write $s=\widetilde s\otimes e_L^{\otimes k}$ locally, where $\widetilde s$ is an $F$-valued differential form. Then define $|s|_0:=|\widetilde s|_{\omega,h_F}$. For $D\subset D_1$, write
\[
 \|s\|_{k\phi_\nu,D}^2=\int_D|s|_0^2{e^{-k\phi_\nu}}\,dV,
 \quad
 \|s\|_{k\phi,D}^2=\int_D|s|_0^2{e^{-k\phi}}\,dV.
\]

Write $\dbar^*_{k\phi}$ and $\dbar^*_{k\phi_\nu}$ for the Hilbert adjoints of the maximal local Dolbeault operators in these weighted spaces.

Choose $D_0\Subset D_1$ and a cut-off function $\chi\in C_c^\infty(D_1)$ with $\chi=1$ near $\overline{D_0}$. Note that
\[
\bar\partial (\chi u)=
\bar\partial\chi\wedge u+\chi\bar\partial u=
\bar\partial\chi\wedge u.
\]
For $a\in\Dom\dbar$, we have
\[
 0=\langle\dbar(\chi a),u\rangle_{k\phi}
 =\langle\dbar a,\chi u\rangle_{k\phi}
+\langle a,(\dbar\chi\wedge)^*u\rangle_{k\phi}.
\]
Then $\chi u\in\Dom\dbar\cap\Dom\dbar^*_{k\phi}$ and $\dbar^*_{k\phi}(\chi u)=-(\dbar\chi\wedge)^*u$. Since $\chi=1$ near $D_0$, $\dbar (\chi u)$ and $\dbar^*_{k\phi}(\chi u)$ vanish on $D_0$. Now we can apply Lemma \ref{lem8} and get $(n,q)$-forms $u_\nu$ compactly supported on $D_1$ such that
\[
 \begin{aligned}
 {e^{-k\phi_\nu/2}}u_\nu\longrightarrow{e^{-k\phi/2}}\chi u, \;
 {e^{-k\phi_\nu/2}}\dbar u_\nu\longrightarrow{e^{-k\phi/2}}\dbar (\chi u),\;
 {e^{-k\phi_\nu/2}}\dbar^*_{k\phi_\nu}u_\nu
     \longrightarrow{e^{-k\phi/2}}\dbar^*_{k\phi}(\chi u)
 \end{aligned}
\]
in $L^2(D_1)$. This means, for example, $\int_{D_1}\left|{e^{-k\phi_\nu/2}}u_\nu-{e^{-k\phi/2}}(\chi u) \right|_0^2\,dV\rightarrow 0$. In particular, we have in $L^2(D_0)$
\begin{equation}\label{eq11}
 \|\dbar u_\nu\|_{k\phi_\nu,D_0}
 +\|\dbar^*_{k\phi_\nu}u_\nu\|_{k\phi_\nu,D_0}\longrightarrow0,\quad
 {e^{-k\phi_\nu/2}}u_\nu\longrightarrow{e^{-k\phi/2}}u.
\end{equation}

Put $\gamma_\nu=*_\omega u_\nu$ with $T_\nu=T_{u_\nu}$ formed using $e^{-k\phi_\nu}$. We deduce from \eqref{eq11} and the Cauchy-Schwarz inequality that
\begin{equation}\label{eq12}
 T_\nu\longrightarrow T_u,\quad
 |u_\nu|_0^2{e^{-k\phi_\nu}}\longrightarrow|u|_0^2{e^{-k\phi}}
\end{equation}
in $L^1(D_0)$, Let $f \in C_c^\infty(D_0)$ be nonnegative, which yields
\begin{equation}\label{eq13}
 \bigl\langle\ddc T_\nu\wedge\omega_{q-1},f\bigr\rangle
 =\int_{D_0}T_\nu\wedge\ddc f\wedge\omega_{q-1}
 \longrightarrow
 \bigl\langle\ddc T_u\wedge\omega_{q-1},f\bigr\rangle.
\end{equation}

Applying Lemma \ref{lem3} to $u_\nu$ with the local weight ${e^{-k\phi_\nu}}$, we write
\begin{equation}\label{eq14}
\begin{aligned}
\int f|\bar\partial\gamma_\nu|_0^2{e^{-k\phi_\nu}}\,dV
\le&\left\langle
\ddc T_\nu\wedge\omega_{q-1},f
\right\rangle+C\int f|u_\nu|_0^2{e^{-k\phi_\nu}}\,dV\\
&+\int f|\bar\partial u_\nu|_0^2{e^{-k\phi_\nu}}\,dV+2\left|
\left\langle\dbar^*_{k\phi_\nu}u_\nu,
\dbar^*_{k\phi_\nu}(fu_\nu)
\right\rangle_{k\phi_\nu}\right|.
\end{aligned}
\end{equation}
We have by \eqref{eq11}
\begin{equation}\label{eq15}
0\leq\int f|u_\nu|_0^2{e^{-k\phi_\nu}}\,dV\to \int f|u|_0^2{e^{-k\phi}}\,dV.
\end{equation}
Note that $\dbar^*_{k\phi_\nu}(fu_\nu)=f\dbar^*_{k\phi_\nu}u_\nu-(\dbar f\wedge)^*u_\nu$. Consequently with \eqref{eq11}, we obtain
\begin{equation}\label{eq16}
 \bigl|\langle\dbar^*_{k\phi_\nu}u_\nu,
\dbar^*_{k\phi_\nu}(fu_\nu)\rangle_{k\phi_\nu}\bigr|
\le\|\dbar^*_{k\phi_\nu}u_\nu\|_{k\phi_\nu,D_0}
 \bigl(\|f\|_\infty\|\dbar^*_{k\phi_\nu}u_\nu\|_{k\phi_\nu,D_0}
+C\|df\|_\infty\|u_\nu\|_{k\phi_\nu,D_0}\bigr)
 \rightarrow0.
\end{equation}

Let $f\equiv 1$ near a compact subset $K\Subset D_0$ in \eqref{eq14}. Then combined with \eqref{eq12}, \eqref{eq13}, \eqref{eq15} and \eqref{eq16}, we have
\begin{equation}\label{eq17}
c\sup_\nu\int_K|\dbar\gamma_\nu|_0^2\,dV\leq\sup_\nu\int_K|\dbar\gamma_\nu|_0^2{e^{-k\phi_\nu}}\,dV<\infty.
\end{equation}

We claim that the weak limit of $\bar\partial\gamma_\nu$ is  $\bar\partial\gamma$. It is easy to see that $u_\nu\to u$ in $L^2(D_0)$ by \eqref{eq11}. Thus $\gamma_\nu\to\gamma$ locally in $L^2$ and $\bar\partial\gamma_\nu\rightarrow\bar\partial\gamma$ in the sense of currents. It follows from \eqref{eq17} that any weak subsequential limit of $\{\bar\partial\gamma_\nu\}$ is $\bar\partial\gamma$, which proves the claim.

For fixed $\mu$, recall that ${e^{-k\phi_\mu}}\le{e^{-k\phi_\nu}}$ for $\nu\ge\mu$. The weak lower semicontinuity of $L^2$ norm \cite[Proposition 3.5]{Brezis} yields
\[
 \int f|\dbar\gamma|_0^2{e^{-k\phi_\mu}}\,dV
 \le\liminf_{\nu\to\infty}
       \int f|\dbar\gamma_\nu|_0^2{e^{-k\phi_\nu}}\,dV.
\]
Letting $\mu\to\infty$ and using the monotone convergence theorem yield
\begin{equation}\label{eq18}
 \int f|\dbar\gamma|_0^2{e^{-k\phi}}\,dV
 \le\liminf_{\nu\to\infty}
\int f|\dbar\gamma_\nu|_0^2{e^{-k\phi_\nu}}\,dV.
\end{equation}
This proves $\bar\partial\gamma\in L^2_{\mathrm{loc}}$. With \eqref{eq12}, \eqref{eq13}, \eqref{eq15}, \eqref{eq16} and \eqref{eq18} in hands, taking limit in \eqref{eq14}, we get
\[
\int f|\bar\partial\gamma|_0^2{e^{-k\phi}}\,dV
\le{}\left\langle
\ddc T_u\wedge\omega_{q-1},f
\right\rangle+C\int f|u|_0^2{e^{-k\phi}}\,dV.
\]
Then $\dbar\gamma\in L^2_{p,1}(X,F\otimes L^k,h^k)$ and \eqref{eq9} follow from a partition of unity. Finally testing \eqref{eq9} with the constant function $1$ yields \eqref{eq10}. The proof is completed.
\end{proof}

\subsection{Local mass estimates}
For a fixed $x\in X$, put $\mu=T_u\wedge\omega_q$ and $\sigma(r)=\mu(B(x,r))$ in this paragraph. To obtain the local mass estimates, one crucial idea is to test the current inequality~\eqref{eq19} with a radial cut-off function to get \eqref{eq24}, in which the interior integral involves the trace measure of $T_u$, the boundary integral involves the rank-one form $d\varrho\wedge d^c\varrho$. This will enable us to establish $\sigma'(r)\ge(2q/r-C)\sigma(r)$ for $\sigma(r)=\mu(B(x,r))$, and hence the mass bound of order $r^{2q}$. The ideas of this subsection is inspired by \cite{Berndtsson}.

\begin{proposition}\label{prop2}
Let $T$ be a positive current of bidegree $(p,p)$ with locally integrable coefficients on $X$. Denote by $\mu=T\wedge\omega_q$ the trace measure of $T$. Suppose that for some $C_0\ge0$,
\begin{equation}\label{eq19}
 \ddc T\wedge\omega_{q-1}\ge-C_0\mu.
\end{equation}
Then there are constants $r_0,C>0$, depending only on $\omega, C_0$, such that
\begin{equation}\label{eq20}
 \mu(B(x,r))\le Cr^{2q}\mu(X),
\end{equation}
where $0<r<r_0$, $B(x,r)$ denotes a ball in a local chart on $X$.
\end{proposition}

\begin{proof}
Set $S:=T\wedge\omega_{q-1}$. Then \eqref{eq19} is equivalent to $\ddc S\ge-C_0\mu$. We have $ S\wedge\omega=q\mu$. If $\alpha$ is a semipositive $(1,1)$-form of rank one satisfying $0\le\alpha\le\omega$, we have
\begin{equation}\label{eq21}
S\wedge\alpha\le\mu.
\end{equation}
To this end, choose a unitary $(1,0)$-coframe $\theta_1,\ldots,\theta_n$ adapted to $\alpha$, put $e_1=i\theta_1\wedge\overline{\theta_1}$ and $\omega_\perp=i\sum_{j=2}^n\theta_j\wedge\overline{\theta_j}$, and write $\omega=e_1+\omega_\perp$ and $\alpha=ae_1$, with $0\le a\le1$. Then
\[
 \omega_q-\alpha\wedge\omega_{q-1}
 =(1-a)e_1\wedge\frac{\omega_\perp^{q-1}}{(q-1)!}
 +\frac{\omega_\perp^q}{q!}
\]
is strongly positive. Then \eqref{eq21} follows by pairing with $T$.

Fix $x\in X$ and choose a local chart with $x=0, \omega=\beta$ at $x$. Let $\varrho=|z|^2$, then $\ddc \varrho=2\beta$ For sufficiently small $r$, there is a constant $C_1>0$, on $B(x,r)$ such that
\begin{equation}\label{eq22}
 (1-C_1r)\omega\le\beta\le(1+C_1r)\omega.
\end{equation}

Let $d\lambda$ be standard Lebesgue measure on $\C^n$ and $dS$ Euclidean surface measure on coordinate spheres. Set $\sigma(r):=\mu(B(x,r))$ and write $d\mu=m(z)\,d\lambda(z)$ in the local chart, where $m$ is locally integrable. We have
\[
\sigma(r)=\int_{|z|<r}m(z)\,d\lambda(z)=
\int_0^r\left(\int_{|z|=s}m(z)\,dS(z)\right)ds.
\]
Thus $\sigma$ is absolutely continuous locally, and $\sigma'(r)=\int_{|z|=r}m(z)\,dS(z)$ for a.e. $r>0$.

Consider the radial function
\[
f_r(z):=(r^2-\varrho(z))_+=
\begin{cases}
r^2-|z|^2,& |z|<r,\\
0,& |z|\ge r.
\end{cases}
\]
Write $\delta$ for the Dirac delta in the scalar radial variable. It is easy to compute that
\begin{equation}\label{eq23}
dd^cf_r=-2\mathbf1_{B(x,r)}\beta+\delta(\varrho-r^2)d\varrho\wedge d^c\varrho.
\end{equation}
By the assumption, \eqref{eq23} and Section \ref{sec8}, we get for a.e. $r>0$
\begin{equation}\label{eq24}
-2\int_{B(x,r)}S\wedge\beta+
\int\delta(\varrho-r^2)S\wedge d\varrho\wedge d^c\varrho=\langle S,dd^cf_r\rangle=\langle dd^cS,f_r\rangle\ge -C_0\int f_r\,d\mu.
\end{equation}
Hence for a.e. $r>0$
\begin{equation}\label{eq25}
2\int_{B(x,r)}S\wedge\beta\le
\int\delta(\varrho-r^2)S\wedge d\varrho\wedge d^c\varrho
+C_0r^2\sigma(r).
\end{equation}

It follows from \eqref{eq22} that
\begin{equation}\label{eq26}
q(1-C_1r)\sigma(r)\le \int_{B(x,r)}S\wedge\beta.
\end{equation}
For the first term in the RHS of \eqref{eq25}, $\frac{1}{2\varrho}d\varrho\wedge d^c\varrho$ is a semipositive $(1,1)$-form of rank one. Note that $0\le\frac{1}{2\varrho}d\varrho\wedge d^c\varrho\le \beta\le(1+C_1r)\omega$. By applying \eqref{eq21} we deduce that on the sphere $\varrho=r^2$
\[
S\wedge d\varrho\wedge d^c\varrho\le2r^2(1+C_1r)\mu.
\]
Therefore by the coarea formula \cite{Hor, PSWW}, we have for a.e. $r>0$
\begin{equation}\label{eq27}
\int\delta(\varrho-r^2)S\wedge d\varrho\wedge d^c\varrho
\le 2r^2(1+C_1r)
\int\delta(\varrho-r^2)\,d\mu=r(1+C_1r)\sigma'(r).
\end{equation}
We infer from \eqref{eq24}, \eqref{eq26} and \eqref{eq27} that
\[
2q(1-C_1r)\sigma(r)\le r(1+C_1r)\sigma'(r)+C_0r^2\sigma(r),
\]
that is,
\[
\sigma'(r)\ge\frac{2q(1-C_1r)-C_0r^2}{r(1+C_1r)}\sigma(r)\ge (\frac{2q}{r}-C_2)\sigma(r),
\]
for some constant $C_2>0$. Then a direct computation shows
\[
\sigma(r)\le e^{C_2(r_0-r)}\left(\frac{r}{r_0}\right)^{2q}\sigma(r_0),
\]
for $0<r<r_0$. We fix $r_0$. Note that $\sigma(r_0)\le\mu(X)$. Finally we get
\[
\mu(B(x,r))=\sigma(r)\le Cr^{2q}\mu(X)
\]
for some constant $C>0$. The proof is completed.
\end{proof}

\begin{corollary}\label{cor2}
Let $1\le q\le n$, $u\in\HH_k^{n,q}(F)$. Then there are constants $C,r_0>0$ depending only on $\omega, h_F$, such that
\begin{equation}
 \int_{B(x,r)}|\gamma|_{h^k}^2\,dV
       =\int_{B(x,r)}|u|_{h^k}^2\,dV
       \le Cr^{2q}\|u\|_{h^k}^2
\end{equation}
for $x\in X,\ 0<r<r_0$. For $A>1$ there are also constants $C_A,r_A>0$ depending only on $A, \omega, h_F$, such that
\begin{equation}\label{eq29}
 \int_{B(x,2Ar)}|\dbar\gamma|_{h^k}^2\,dV
       \le C_A r^{2q-2}\|u\|_{h^k}^2
\end{equation}
for $x\in X,\ 0<r<r_A<1$.
\end{corollary}

\begin{proof}
Write $\mu=T_u\wedge\omega_q=|u|_{h^k}^2dV$. By Proposition \ref{prop1}, we have
\[
dd^cT_u\wedge\omega_{q-1}\ge-C\mu. 
\]
Hence we can apply Proposition \ref{prop2} and get
\[
\int_{B(x,r)}|\gamma|_{h^k}^2\,dV=\int_{B(x,r)}|u|_{h^k}^2\,dV=\mu(B(x,r))\le Cr^{2q}\mu(X)= Cr^{2q}\|u\|_{h^k}^{2}
\]
for $0<r<r_0$.

Fix $A>1$, choose $r$ sufficiently small so that $B(x,3Ar)$ is contained in the local chart and $3Ar<r_0$. Let
\[
 \eta(z)=\chi\!\left(\frac{|z|^2}{A^2r^2}\right),
\]
where $\chi\in C^\infty(\R,[0,1])$ equals $1$ on $(-\infty,4]$ and $0$ on $[8,\infty)$. Thus $\eta=1$ on $B(x,2Ar)$, $\supp\eta\Subset B(x,3Ar)$ and
\begin{equation}\label{eq30}
 -\widetilde C_Ar^{-2}\omega\le \ddc \eta\le \widetilde C_Ar^{-2}\omega
\end{equation}
for some constant $\widetilde C_A>0$. By using \eqref{eq9} to $\eta$, we have
\begin{equation}\label{eq31}
\begin{split}
 \int_{B(x,2Ar)}|\dbar\gamma|_{h^k}^2\,dV\le \int_{B(x,2Ar)}|\dbar\gamma|_{h^k}^2\eta\,dV&
 \le\int T_u\wedge\omega_{q-1}\wedge \ddc \eta+C\int\eta\,d\mu\\
 &\le\int T_u\wedge\omega_{q-1}\wedge \ddc \eta+C\mu(B(x,3Ar)).
 \end{split}
\end{equation}
It follows from \eqref{eq30} that
\begin{equation}\label{eq32}
 \int T_u\wedge\omega_{q-1}\wedge \ddc \eta\le q\widetilde C_Ar^{-2}\mu(B(x,3Ar)).
\end{equation}
With \eqref{eq20}, \eqref{eq31} and \eqref{eq32} in hands, we conclude that
\[
 \int_{B(x,2Ar)}|\dbar\gamma|_{h^k}^2\,dV
 \le \widetilde C_A(r^{-2}+1)\int_{B(x,3Ar)}|u|_{h^k}^2\,dV
 \le \frac{1}{2}C_A(r^{-2}+1)r^{2q}\|u\|_{h^k}^2.
\]
Taking $r<r_A<1$ proves \eqref{eq29}.
\end{proof}

\section{Gaussian localization and spectral gap in singular metrics}\label{sec4}

Section~\ref{sec3} bounds the mass of a weighted harmonic form $u$ and $\dbar\gamma$ on small balls. To turn these estimates into a cohomology dimension bound for our purpose, we need a finite-dimensional space $\operatorname{Ran}\Pi_{k,\ell}$ (cf. Lemma \ref{lem7}) into which we can project the localized sections of $\gamma$ with a bounded norm below in Section~\ref{sec5}. To this end, we construct in Lemma \ref{lem4} Hermitian line bundles with which the space is defined and Gaussian sections with estimates to control spectral projections. Lemma~\ref{lem5} provides a uniform Bochner type estimate, with which we obtain a spectral gap and bound for the dimension of the lower energy spectral subspace in Lemma~\ref{lem6}. Lemma~\ref{lem7} bounds the rank of the projection $\Pi_{k,\ell}$ by $C_Mk^n$ and norm of complement of $\Pi_{k,\ell}$. Such bounds will be used in the proof of Proposition~\ref{prop3}.

Recall that $h=h_0e^{-U},U\le0$. By the semipositivity of $h$, there exists a uniform constant $C_0>0$ such that
\begin{equation}\label{eq33}
-C_0\omega\leq \Theta_{h_0}(L)\leq C_0\omega, \quad  \ddc U\ge -\Theta_{h_0}(L)\ge -C_0\omega.   
\end{equation}

\subsection{Integral curvature approximation and local Gaussian sections}

The following lemma combines Laeng's approximation of a real cohomology class by integral curvature classes \cite[Part II, Theorem 4.3 and Section 4.3]{LaengThesis} with Popovici's construction of local approximately holomorphic Gaussian sections \cite[Section~3.1]{PopoviciTranscendental}. This result suffices both for the following Bochner type estimate and the control of the norm in \eqref{eq62}.

\begin{lemma}\label{lem4}
For every positive integer $\ell$, there exist a smooth line bundle $H_\ell$ with Hermitian metric $h_{H_\ell}$ and a Hermitian connection $D_\ell$ such that
\begin{equation}\label{eq34}
 iR^{D_\ell}=\ell\omega+e_\ell,
\end{equation}
where $e_\ell$ is a closed $2$-form and $\|e_\ell\|_{\mathcal C^m}\le C_m$ for every integer $m\ge0$, with $C_m$ independent of $\ell$. Moreover, for every $x\in X$, there is a nowhere vanishing smooth section $\sigma_{\ell,x}$ of $H_\ell$ defined on a ball with a uniform radius independent of $\ell, x$ satisfying
\begin{equation}\label{eq35}
 e^{-C\ell d_\omega(x,y)^2}
 \le |\sigma_{\ell,x}(y)|_{h_{H_\ell}}\le e^{-c\ell d_\omega(x,y)^2}, \quad
 |D_\ell''\sigma_{\ell,x}(y)|_{h_{H_\ell}}
 \le C d_\omega(x,y)|\sigma_{\ell,x}(y)|_{h_{H_\ell}},
\end{equation}
where the constants $c,C>0$ are independent of $\ell, x$ and $d_\omega$ denote the distance induced by $\omega$.
\end{lemma}

\begin{proof}
Choose a basis $\gamma_1,\ldots,\gamma_b$ of the lattice $\operatorname{Im}(H^2(X,\Z)\to H^2(X,\R))$ and write $[\omega/(2\pi)]=\sum_{j=1}^b a_j\gamma_j$. Let $\eta_j$ be a smooth real closed form which represents $2\pi\gamma_j$. Choose integers $m_{\ell j}$ satisfying $|m_{\ell j}-\ell a_j|\le1/2$ and set $e_\ell:=\sum_{j=1}^b(m_{\ell j}-\ell a_j)\eta_j$. Then
\begin{equation}\label{eq36}
 \|e_\ell\|_{C^m}\le\frac12\sum_{j=1}^b\|\eta_j\|_{C^m},
 \quad
 [\ell\omega+e_\ell]=2\pi\sum_{j=1}^b m_{\ell j}\gamma_j.
\end{equation}
Choose a smooth Hermitian line bundle $H_\ell$ with the integral first Chern class $[\ell\omega+e_\ell]/(2\pi)$ and a corresponding Hermitian connection $D_\ell^0$. Since $\ell\omega+e_\ell-iR^{D_\ell^0}$ is exact, there exists a real $1$-form $\theta_\ell$ such that
\[
 d\theta_\ell=\ell\omega+e_\ell-iR^{D_\ell^0}.
\]
The Hermitian connection $D_\ell=D_\ell^0-i\theta_\ell$ satisfies \eqref{eq34}.

Fix $x\in X$ and use the normalized holomorphic coordinates $z=z_x$ from Lemma~\ref{lem9} in the appendix. On a coordinate ball of uniformly small radius, choose a local potential $\varphi_x$ satisfying
\begin{equation}\label{eq37}
\ddc\varphi_x=\omega,\quad 
c|z|^2\le\varphi_x(z)\le C|z|^2,
\end{equation}
for some constants $c, C>0$. Consider the Hermitian connection $D_{\ell,x}=d+\frac{\ell}{2} (\bar\partial\varphi_x-\partial\varphi_x)$ with $iR^{D_{\ell,x}}=\ell\omega$. The Poincaré homotopy formula provides locally an imaginary $1$-form $a_{\ell,x}$ satisfying $da_{\ell,x}=-i e_\ell$, and let $\iota_z$ denote contraction with the radial vector determined by $z$. Then $a_{\ell,x}(z)=-i\int_0^1t\,\iota_{z}e_\ell(tz)\,dt$. We see that $|a_{\ell,x}(z)|_{\omega}\le C|z|$ by the inequality in \eqref{eq36}.

Note that $D_\ell$ and $D_{\ell,x}+a_{\ell,x}$ have the same curvature. Then $D_\ell-D_{\ell,x}-a_{\ell,x}$ is a closed imaginary $1$-form and hence equals locally $i\,df$ for a real smooth function $f$. Scaling the frame multiplied by $e^{-if}$, we have
\[
 D_\ell=d+\frac{\ell}{2}(\dbar\varphi_x-\partial\varphi_x)
+a_{\ell,x}.
\]
Let $\tau_{\ell,x}$ be a local unitary frame of $H_\ell$. Now define $\sigma_{\ell,x}(z)=e^{-\ell\varphi_x(z)/2}\tau_{\ell,x}$, Then \eqref{eq37} yields
\[
e^{-C\ell|z|^2}\le|\sigma_{\ell,x}(z)|_{h_{H_\ell}}
\le e^{-c\ell|z|^2}.
\]
By Lemma~\ref{lem9}\textup{(i)}, $|z_x(y)|\approx d_\omega(x,y)$ on the chosen coordinate ball, this proves the first inequality in \eqref{eq35}. Moreover,
\[
 D_\ell''\sigma_{\ell,x}
 =\left(\dbar+\frac{\ell}{2}\dbar\varphi_x
+a_{\ell,x}^{0,1}\right)e^{-\ell\varphi_x/2}
 =a_{\ell,x}^{0,1}\sigma_{\ell,x}.
\]
Since $a_{\ell,x}(z)$ is bounded by $C|z|$, we finally prove \eqref{eq35}.
\end{proof}

\begin{remark}
If $[\omega/(2\pi)]$ is integral, there exists a holomorphic line bundle with a smooth Hermitian metric whose Chern curvature is $\omega$. On a general K\"ahler manifold $[\omega/(2\pi)]$ needs not be integral. We overcome this obstruction by using smooth Hermitian line bundles with Hermitian connections that are not necessarily holomorphic in Lemma \ref{lem4}. The error term $e_\ell$ remains uniformly bounded, even though the curvature term grows with order $\ell$. This yields the Bochner type estimate used below to establish the spectral gap, even when the regularized metrics on $L$ are not semipositive.
\end{remark}

\subsection{A Bochner type estimate for non-integrable connections}

Let $G$ be a holomorphic vector bundle endowed with a smooth Hermitian metric $h_G$. Here $\dbar$ acts on the holomorphic factor $G\otimes L^k$, and $D_\ell^{\prime\prime}$ acts on $H_\ell$; their sum is understood by the tensor-product Leibniz rule. Let $H_{k,\ell}=G\otimes L^k\otimes H_\ell$ and, for $0\le j\le n$, set
\[
A_j=\dbar+D_\ell'':C^\infty(X,\Lambda^{0,j}T^*X\otimes H_{k,\ell})
\rightarrow
C^\infty(X,\Lambda^{0,j+1}T^*X\otimes H_{k,\ell}).
\]
Note that $A=\{A_j\}$ do not form a complex since $A_{j+1}A_j\neq 0$ in general. For the approximations $U_\nu$ in Lemma \ref{lem1} and $0\le s\le1$, equip $H_{k,\ell}$ with the metric
\begin{equation}\label{eq38}
 h_{s,\nu}=h_G \otimes h_0^k \otimes h_{H_\ell}e^{-skU_\nu}.
\end{equation}
Let $\|\cdot\|_{s,\nu}$ and $A_{j}^{*s,\nu}$ denote the corresponding $L^2$ norm and Hilbert adjoint with respect to $h_{s,\nu}$. Let $D_0$ be the Hermitian product connection for $h_{0,\nu}$ with $(D_0)^{\prime\prime}=A$. The connection $D_{s,\nu}$, which is compatible with $h_{s,\nu}$ and $A$, satisfies then
\[
D_{s,\nu}''=A,\quad D_{s,\nu}'=D_0'-sk\,\partial U_\nu.
\]
Let $M_0=4C_0$, where $C_0$ is the constant in \eqref{eq33}. In the following lemma, we establish a Bochner type estimate with the spectrum bound of $A_0^{*s,\nu}A_0$.
\begin{lemma}\label{lem5}
There exist positive constants $c_0, k_0, R_0$ such that, for every integer $M\ge M_0$, $k\ge k_0$, and $\ell=Mk$, we have
\begin{equation}\label{eq39}
 \|A_jv\|_{s,\nu}^2+\|A_{j-1}^{*s,\nu}v\|_{s,\nu}^2
 \ge c_0\ell\|v\|_{s,\nu}^2
\end{equation}
for $(0,j)$-form $v$ with values in $H_{k,\ell}$, $1\le j\le n,\ 0\le s\le 1$. Moreover, for every smooth section $v$ of $H_{k,\ell}$,
\begin{equation}\label{eq40}
 \|A_1A_0v\|_{s,\nu}\le R_0\|v\|_{s,\nu}.
\end{equation}
All constants are independent of $s,\nu,k$ and $M$.
\end{lemma}

\begin{proof}
For $(0,j)$-forms, we use the canonical isometric identification
\[
\Lambda^{0,j}T^*X\otimes H_{k,\ell}\cong
\Lambda^{n,j}T^*X\otimes H_{k,\ell}\otimes K_X^{-1}.
\]
We estimate the curvature of $H_{k,\ell}\otimes K_X^{-1}$ from below. Note that the curvature of $G$ and $K_X$ are smooth. We only consider the $(1,1)$ part of the curvature contribution from $L^k\otimes H_\ell$, which is
\[
k\Theta_{h_0}(L)+\ell\omega+e_\ell^{1,1}+sk\,dd^cU_\nu
\]
by \eqref{eq34} and \eqref{eq38}. It follows from \eqref{eq33} and Lemma \ref{lem1} that
\[
\begin{split}
k\Theta_{h_0}(L)+\ell\omega+sk\,dd^cU_\nu
&\ge-kC_0\omega+\ell\omega-kC_0\omega=(\ell-2C_0k)\omega\\
&=(M-2C_0)k\omega\ge \frac{\ell}{2}\omega,
\end{split}
\]
for $M\ge M_0\ge 4C_0$. Since $e_\ell^{1,1}$ is bounded by Lemma \ref{lem4}, we finally get
\begin{equation}\label{eq41}
iR^{D_{s,\nu}}\ge c_0\ell\omega.
\end{equation}
We write $\Delta':=D_{s,\nu}'(D'_{s,\nu})^*+(D'_{s,\nu})^*D'_{s,\nu}$ and similarly for $\Delta''$. The Bochner-Kodaira-Nakano identity and \eqref{eq41} imply that
\begin{equation}
\begin{split}
\|D''_{s,\nu}v\|_{s,\nu}^2+\|(D''_{s,\nu})^*v\|_{s,\nu}^2&=\langle\Delta''v,v\rangle_{s,\nu}=\langle\Delta'v,v\rangle_{s,\nu}+\left\langle
[iR^{D_{s,\nu}}\wedge,\Lambda]v,v
\right\rangle_{s,\nu} \\
&\ge c_0\ell\left\langle[\omega\wedge,\Lambda]v,v\right\rangle_{s,\nu}=
c_0\ell\,j\,\|v\|_{s,\nu}^2\ge c_0\ell\|v\|_{s,\nu}^2,
\end{split}
\end{equation}
Then \eqref{eq39} follows for $A_j=D''_{s,\nu}$ on $(0,j)$-forms.

Since $G$ and $L$ are holomorphic Hermitian bundles, it is easy to see that $A_1A_0v=(R^{D_\ell})^{0,2}\wedge v$. Applying Lemma \ref{lem4} again, we get $(R^{D_\ell})^{0,2}=-i(e_\ell)^{0,2}$ and $|A_1A_0v|_{s,\nu}^2\le C|v|_{s,\nu}^2$. Hence $\|A_1A_0v\|_{s,\nu}\le R_0\|v\|_{s,\nu}$. The proof is completed.
\end{proof}

\subsection{Spectral gap and rank bound for Hermitian line bundles}
Set $T_{s,\nu}:=A_0^{*s,\nu}A_0$, which is a nonnegative elliptic self-adjoint operator with compact resolvent. In Lemma~\ref{lem6}, we apply Lemma \ref{lem5} to get a spectral gap inspired by the argument in \cite[Corollary~2.2]{PopoviciTranscendental}. The spectral gap allows us to estimate the dimension of the lower energy spectral subspace as shown in Lemma \ref{lem7}.

We choose $c_0, R_0, M, k, \nu$ as in Lemma \ref{lem5}.
\begin{lemma}\label{lem6}
For $0\le s\le1$, $\ell=Mk$ and $k$ large, we have the spectral gap
\begin{equation}\label{eq43}
 \Spec T_{s,\nu}
 \subset[0,\frac{2R_0^2}{c_0\ell}]
\cup[c_0\ell/2,\infty).
\end{equation}
Moreover, the rank $N_0(k,\ell)$ of $\one_{[0,c_0\ell/4]}(T_{s,\nu})$ is independent of $s,\nu$ and satisfies
\begin{equation}\label{eq44}
 N_0(k,\ell)\le C(k+\ell)^n.
\end{equation}
\end{lemma}

\begin{proof}
Let $T_{s,\nu}v=\lambda v$, $v\ne0$. A direct computation and \eqref{eq40} yield
\[
 \|A_0v\|_{s,\nu}^2=\lambda\|v\|_{s,\nu}^2,\quad
 \|A_0^{*s,\nu}A_0v\|_{s,\nu}^2=\lambda^2\|v\|_{s,\nu}^2,\quad
 \|A_1A_0v\|_{s,\nu}^2\le R_0^2\|v\|_{s,\nu}^2,
\]
where the norm $\|\cdot\|_{s,\nu}$ is with respect to $\omega, h_{s,\nu}$. Then combined with\eqref{eq39} we obtain
\[
 c_0\ell\lambda\|v\|_{s,\nu}^2
 \le\|A_1A_0v\|_{s,\nu}^2+\|A_0^{*s,\nu}A_0v\|_{s,\nu}^2
 \le(\lambda^2+R_0^2)\|v\|_{s,\nu}^2.
\]
Thus $\lambda(c_0\ell-\lambda)\le R_0^2$. If $\lambda<c_0\ell/2$, then $\lambda\le2R_0^2/(c_0\ell)$, proving \eqref{eq43}.

It is known that
\begin{equation}\label{eq45}
\|(T_{s,\nu}+1)^{-1}-(T_{s',\nu}+1)^{-1}\|_{\mathrm{op},h_{0,\nu}}\rightarrow0
\end{equation}
as $s'\to s$ (cf. \cite[Chapter~VI, Theorem~3.6]{Kato}). For sufficiently large $\ell$, $2R_0^2/(c_0\ell)<c_0\ell/4$. Thus $c_0\ell/4$ belongs to the spectral gap. Then we can choose a closed contour $\Gamma$ surrounding $[0,\frac{2R_0^2}{c_0\ell}]$ and separating it from $[c_0\ell/2,\infty)$. Therefore the spectral projection is
\[
P_{s,\nu}=\frac{1}{2\pi i}\int_\Gamma(z-T_{s,\nu})^{-1}\,dz,
\]
where we write $P_{s,\nu}=\mathbf1_{[0,c_0\ell/4]}(T_{s,\nu})$ for short. Moreover, we have by \eqref{eq45}
\[
\|P_{s,\nu}-P_{s',\nu}\|_{\mathrm{op},h_{0,\nu}}\rightarrow0.
\]
It is well-known that if two orthogonal projections $P,Q$ satisfy $\|P-Q\|_{\mathrm{op}}<1$, then their ranges have the same dimension. Thus $N_0(k,\ell)$ is constant in $s\in [0,1]$. At $s=0$, there is no $U_\nu$ involved in the metric. Hence $N_0(k,\ell)$ is also independent of $\nu$.

Now it suffices to show \eqref{eq44} only when $s=0$. Write $T_0=T_{0,\nu}$ for short. Let $\nabla=D_{0,\nu}$ be the Hermitian connection on $H_{k,\ell}$. Note that $\nabla^*\nabla=\Delta'+\Delta''$ and $\Delta''=\Delta'-i\Lambda R^\nabla$. Hence we have the Weitzenb\"ock formula
\begin{equation}\label{eq46}
 T_0=\frac12\nabla^*\nabla-\frac{i}{2}\Lambda R^\nabla.
\end{equation}
Recall that the curvature of $H_{k,\ell}$ is $R^{H_{k,\ell}}=R^G+kR^L+R^{H_\ell}$. Then we have
\begin{equation}\label{eq47}
\|\Lambda R^\nabla\|_\infty\le C(k+\ell)
\end{equation}
by using \eqref{eq34}.

We now compare the traces of the heat kernels of $T_0$ and $\nabla^*\nabla/2$. Denote by $\{\lambda_j\}$ and $\{\mu_j\}$ the eigenvalues of $T_0$ and $\nabla^*\nabla/2$, respectively, listed in nondecreasing order with multiplicity. Recall that the min-max principle is
\[
\lambda_j=\inf_{\dim V=j}\sup_{0\ne v\in V}
\frac{\langle T_0v,v\rangle_{h_{0,\nu}}}{\|v\|_{h_{0,\nu}}^2}.
\]
It follows from \eqref{eq46} and \eqref{eq47} that
\[
\langle T_0v,v\rangle_{h_{0,\nu}}\ge
\langle\tfrac12\nabla^*\nabla v,v\rangle_{h_{0,\nu}}-C(k+\ell)\|v\|_{h_{0,\nu}}^2.
\]
Hence by applying the above min-max principle we deduce that
\[
\lambda_j\ge\mu_j-C(k+\ell).
\]
Then for $t>0$,
\begin{equation}\label{eq48}
\operatorname{Tr}(e^{-tT_0})\le e^{Ct(k+\ell)}
\operatorname{Tr}e^{-\frac{t}{2}\nabla^*\nabla}.
\end{equation}
It follows from the semigroup domination theorem for Hermitian connections \cite[Theorem 2.13]{Guneysu} that
\[
\operatorname{Tr}e^{-t\nabla^*\nabla}(x,x)\le (\operatorname{rank}G) e^{-t\Delta}(x,x),
\]
where $\Delta=d^*d+dd^*$. Hence
\begin{equation}\label{eq49}
\begin{split}
\operatorname{Tr}e^{-\frac{t}{2}\nabla^*\nabla}&=\int_X\operatorname{Tr}e^{-t\nabla^*\nabla}(\frac{t}{2},x,x)dV(x)\\
&\le \operatorname{rank}G \int_Xe^{-t\Delta}(\frac{t}{2},x,x)dV(x)=(\operatorname{rank}G) \operatorname{Tr}e^{-\frac{t}{2}\Delta}.
\end{split}
\end{equation}
The small time asymptotics of heat kernels \cite[Theorem 2.41]{BerlineGetzlerVergne} yield for $0<t\le 1$
\begin{equation}\label{eq50}
\operatorname{Tr}(e^{-\frac{t}{2}\Delta})\le Ct^{-n}.
\end{equation}
With \eqref{eq48}, \eqref{eq49} and \eqref{eq50} in hands, we get for $0<t\le 1$
\begin{equation}\label{eq51}
\operatorname{Tr}(e^{-tT_0})\le Ce^{Ct(k+\ell)}t^{-n}.
\end{equation}

In $\operatorname{Ran}\mathbf1_{[0,c_0\ell/4]}(T_0)$, the eigenvalues satisfy $0\le\lambda_j\le\frac{c_0\ell}{4}$. Then
\begin{equation}\label{eq52}
\operatorname{Tr}(e^{-tT_0})=
\sum_{j=1}^\infty e^{-t\lambda_j}\ge
\sum_{j=1}^{N_0(k,\ell)}e^{-t\lambda_j}\ge
N_0(k,\ell)e^{-tc_0\ell/4}.
\end{equation}
We finally obtain from \eqref{eq51} and \eqref{eq52} that
\[
N_0(k,\ell)\le Ce^{C_1t(k+\ell)}t^{-n}.
\]
Then \eqref{eq44} follows immediately by taking $t=(k+\ell)^{-1}\le 1$. The proof is completed.
\end{proof}

\begin{remark}
Note that if $H_\ell$ are holomorphic line bundles, $A_1A_0=0$ and $R_0=0$. Then \eqref{eq43} is reduced to the classical spectral gap condition for positive line bundles. The exponent $n$ in \eqref{eq44} comes from the small time asymptotic expansion of heat kernels, which is fundamental in spectral geometry and index theory. Recently Hsiao, Huang and the second named author \cite{HHS} established the small time asymptotics of heat kernel when the Laplacian is not hypoelliptic. They also deal with a variant of spectral gap condition needed for asymptotics of analytic torsions associated with high powers of a line bundle.
\end{remark}

Now we come back to the singular metric $\widetilde h=h_G\otimes h_0^k\otimes h_{H_\ell}e^{-kU}, U\le0$ of $H_{k,\ell}$. Define
\[
A_0:L^2(X,H_{k,\ell},\widetilde h)\rightarrow
L^2_{0,1}(X,H_{k,\ell},\widetilde h).
\]
We know that $A_0$ is closed and densely defined by employing the same argument in Lemma \ref{lem2}. Write $A_{0,\widetilde h}^*$ for the Hilbert adjoint of $A_0$ in the metric $\widetilde h$. Let $\Pi_{k,\ell}=\mathbf1_{[0,c_0\ell/4]}(A_{0,\widetilde h}^*A_0)$.

\begin{lemma}\label{lem7}
Fix an integer $M\ge M_0$, let $\ell=Mk$. Then we have for sufficiently large $k$
\begin{equation}
 \rank\Pi_{k,\ell}
 \le N_0(k,\ell)\le C_Mk^n,
\end{equation}
and
\begin{equation}\label{eq54}
 \|(I-\Pi_{k,\ell})v\|_{\widetilde h}
 \le\frac{2}{\sqrt{c_0\ell}}\|A_0v\|_{\widetilde h},
\end{equation}
where $v\in\operatorname{Dom}A_0$.
\end{lemma}

\begin{proof}
Let $E_\lambda=\one_{[0,\lambda]}(A_{0,\widetilde h}^*A_0)$ be its spectral family. The spectral theorem shows for $v\in\operatorname{Dom}A_0$
\[
\begin{split}
\|(I-\Pi_{k,\ell})v\|_{\widetilde h}^2&=\int_{\frac{c_0\ell}{4}}^\infty
d\langle E_\lambda v,v\rangle_{\widetilde h}\le \frac{4}{c_0\ell}
\int_{\frac{c_0\ell}{4}}^\infty
\lambda\,d\langle E_\lambda v,v\rangle_{\widetilde h}\\
&\le \frac{4}{c_0\ell}
\int_{0}^\infty
\lambda\,d\langle E_\lambda v,v\rangle_{\widetilde h}=\frac{4}{c_0\ell}\|A_0v\|_{\widetilde h}^2,
\end{split}
\]
which proves \eqref{eq54}.

Let $W\subset\Ran\Pi_{k,\ell}$ be a finite dimensional subspace. Let $d=\dim W$. Because the spectral support of $v\in W$ is contained in $[0,c_0\ell/4]$, we have for $v\in W$
\begin{equation}\label{eq55}
\|A_0v\|_{\widetilde h}^2\le (c_0\ell/4)\|v\|_{\widetilde h}^2.
\end{equation}
Choose an orthonormal basis $v_1,\ldots,v_{d}$ with respect to $\widetilde h$. For $\widetilde h_\nu=h_G\otimes h_0^k\otimes h_{H_\ell}e^{-kU_\nu}$, Lemma \ref{lem1} and the dominated convergence theorem yield $\|v\|_{\widetilde h_\nu}^2\rightarrow\|v\|_{\widetilde h}^2$.

Define the Hermitian matrices as follows
\[
\mathscr G_\nu=\left(\langle v_i,v_j\rangle_{\widetilde h_\nu}
\right)_{i,j=1}^d, \quad \mathscr H_\nu
=\left(\langle A_0v_i,A_0v_j\rangle_{\widetilde h_\nu}
\right)_{i,j=1}^d, \quad \mathscr H=
\left(\langle A_0v_i,A_0v_j\rangle_{\widetilde h}\right)_{i,j=1}^d.
\]
Then the dominated convergence theorem again shows $\mathscr G_\nu\rightarrow I$ and $\mathscr H_\nu\rightarrow\mathscr H$. Moreover, \eqref{eq55} implies that $\mathscr H\le \frac{c_0\ell}{4} I$. Choose $\nu$ large enough that Lemma \ref{lem6} applies and
\[
 \mathscr G_\nu\ge\frac34 I,\quad
 \mathscr H_\nu\le\frac{5}{16}c_0\ell I.
\]
Then for every nonzero $v\in W$, the following holds
\begin{equation}\label{eq56}
 \frac{\|A_0v\|_{\widetilde h_\nu}^2}{\|v\|_{\widetilde h_\nu}^2}
 <\frac{c_0\ell}{2}.
\end{equation}

For the smooth operator $T_{1,\nu}=A_0^{*1,\nu}A_0$ in Lemma \ref{lem6}, we claim that the restriction of $\mathbf1_{[0,c_0\ell/4]}(T_{1,\nu})$ to $W$ is injective. To this end, if $\mathbf1_{[0,c_0\ell/4]}(T_{1,\nu}) v=0$ for a nonzero $v\in W$, then \eqref{eq43} forces $v$ to be in the spectral subspace with spectral support in $[c_0\ell/2,\infty)$. Then the spectral theorem implies
\[
\|A_0v\|_{\widetilde h_\nu}^2\ge\frac{c_0\ell}{2}\|v\|_{\widetilde h_\nu}^2,
\]
which is a contradiction with \eqref{eq56}. By using Lemma \ref{lem6} again, we get $\dim W\le N_0(k,\ell)$. Since $W$ is an arbitrary subspace, we have $\rank\Pi_{k,\ell}\le N_0(k,\ell)\le C_Mk^n$. This completes the proof.
\end{proof}

\section{Proof of the main theorems}\label{sec5}
In this section, we are devoted to proving the main results by developing a spectral localization approach. Proposition~\ref{prop3} is the main step in the proof of Theorem~\ref{thm1}. Fix $1\le q\le n$. Recall that $u\in \HH_k^{n,q}(F)$ and $\gamma=*_\omega u$. For $r=k^{-1/2}$, choose the covering $B_i=B(x_i,r)$ of Lemma \ref{lem9}. Fix $A>2$. Let $D_i=B(x_i,2Ar)$ contained in a local chart. Choose cut-off functions $\chi_i$ satisfying
\[
 0\le\chi_i\le1,\quad
 \chi_i=1\ \hbox{on }B(x_i,Ar),\quad
 \operatorname{supp}\chi_i\Subset D_i,\quad
 |d\chi_i|_{\omega}\le C_A r^{-1}.
\]
There is an overlap bound $N_A$, independent of $k$, such that $\sum_i\one_{D_i}\le N_A$ on $X$ by Lemma \ref{lem9}.

Let $\ell=Mk$ for a fixed sufficiently large integer $M$. Write $\sigma_i=\sigma_{\ell,x_i}$ for the sections in Lemma \ref{lem4}, and $\Pi=\Pi_{k,\ell}$ for the spectral projection in Lemma \ref{lem7}. Let $ \mathcal E_i u=\chi_i\sigma_i\gamma$. Define a local projection by
\begin{equation}
 \mathcal T u=(\Pi\mathcal E_i u)_i:
 \HH_k^{n,q}(F)\rightarrow\bigoplus_i\operatorname{Ran}\Pi.
\end{equation}

\begin{proposition}\label{prop3}
There are $a>0$, $C>0$ independent of $k, u$ such that, for sufficiently large $k$, every $u\in\HH_k^{n,q}(F)$, and every orthonormal family $u_1,\ldots,u_{d}\in\HH_k^{n,q}(F)$, we have
\begin{equation}\label{eq58}
 \|\mathcal T u\|_{\widetilde h}^2\ge a\|u\|_{h^k}^2,\quad
 \sum_{\alpha=1}^{d}\|\mathcal T u_\alpha\|_{\widetilde h}^2\le Ck^{n-q},
\end{equation}
for some suitable $A, M$.
\end{proposition}

\begin{proof}
We first estimate the norm of $A_0\mathcal E_i u$. Consider
\begin{equation}
A_0(\mathcal E_i u)=A_0(\chi_i\sigma_i\gamma)
 =\bar\partial\chi_i\otimes\sigma_i\gamma
   +\chi_i\sigma_i\bar\partial\gamma
   +\chi_i\gamma D_\ell''\sigma_i.
\end{equation}
Then we have
\[
\|A_0\mathcal E_i u\|_{\widetilde h}^2 \leq C\big(\|d\chi_i\otimes\sigma_i\gamma\|_{\widetilde h}^2
+\|\chi_i\sigma_i\bar\partial\gamma\|_{\widetilde h}^2
+\|\chi_i\gamma D_\ell''\sigma_i\|_{\widetilde h}^2\big).
\]
Since $|d\chi_i|_{\omega}\le C_A r^{-1}$, we have
\[
\|d\chi_i\otimes\sigma_i\gamma\|_{\widetilde h}^2
\leq C_A r^{-2}
\int_{U_i}
|\sigma_i|_{h_{H_\ell}}^2|\gamma|_{h^k}^2\,dV,
\]
where $U_i=\operatorname{supp}d\chi_i$. For $x\in U_i$, $d_\omega(x,x_i)\ge Ar$, then by Lemma \ref{lem4}
\[
|\sigma_i(x)|_{h_{H_\ell}}^2=|\sigma_{\ell,x_i}(x)|_{h_{H_\ell}}^2\le e^{-2c\ell d_\omega(x,x_i)^2}\le e^{-2c\ell A^2r^2}\le e^{-cMA^2}.
\]
Hence we have
\begin{equation}\label{eq60}
\begin{split}
\sum_i\|d\chi_i\otimes\sigma_i\gamma\|_{\widetilde h}^2
&\leq C_Ar^{-2}e^{-cMA^2}
\sum_i\int_{D_i}|\gamma|_{h^k}^2\,dV\\
&\le N_AC_A r^{-2}e^{-cMA^2}\|u\|_{h^k}^2,
\end{split}
\end{equation}
where $N_A$ is the overlap bound of the covering $D_i$. For the second term, $\sigma_i$ is bounded. It follows from Proposition \ref{prop1} that
\begin{equation}\label{eq61}
\sum_i\|\chi_i\sigma_i\bar\partial\gamma\|_{\widetilde h}^2
\leq C\sum_i\int_{D_i}
|\bar\partial\gamma|_{h^k}^2\,dV\leq N_AC_A\|u\|_{h^k}^2.
\end{equation}
Lemma \ref{lem4} implies also $|D_\ell''\sigma_i|_{h_{H_\ell}}\leq C_A r|\sigma_i|_{h_{H_\ell}}$ on $D_i$. Consequently,
\begin{equation}\label{eq62}
\sum_i\|\chi_i\gamma D_\ell''\sigma_i\|_{\widetilde h}^2\leq C_A r^2
\sum_i\int_{D_i}|\gamma|_{h^k}^2\,dV\leq N_AC_A r^2\|u\|_{h^k}^2.
\end{equation}
Combining\eqref{eq60}, \eqref{eq61} and \eqref{eq62}, we obtain
\begin{equation}\label{eq63}
\sum_i\|A_0\mathcal E_i u\|_{\widetilde h}^2\leq
2N_AC_A\left(1+r^{-2}e^{-cMA^2}\right)\|u\|_{h^k}^2.
\end{equation}
In particular, $\mathcal E_iu\in \operatorname{Dom}A_0$. It follows from \eqref{eq63} and Lemma \ref{lem7} that
\begin{equation}\label{eq64}
\sum_i\|(I-\Pi)\mathcal E_i u\|_{\widetilde h}^2
\leq \frac{4}{c_0\ell}\sum_i\|A_0\mathcal E_iu\|_{\widetilde h}^2\leq
\frac{\widetilde C_A}{Mk}
\left(1+r^{-2}e^{-cMA^2}\right)\|u\|_{h^k}^2
\end{equation}
for some constant $\widetilde C_A>0$. Recall that $X=\cup_iB(x_i,r)$. For every $x\in X$, choose $x_i$ with $x\in B(x_i,r)$. By applying Lemma \ref{lem4} again, we have
\[
|\sigma_i(x)|_{h_{H_\ell}}=|\sigma_{\ell,x_i}(x)|_{h_{H_\ell}}\ge e^{-C\ell d_\omega(x,x_i)^2}\ge e^{-CkMr^2}=e^{-CM}.
\]
Then $\sum_i|\chi_i(x)\sigma_i(x)|_{h_{H_\ell}}^2\ge e^{-2CM}$. Consequently,
\begin{equation}\label{eq65}
\sum_i\|\mathcal E_i u\|_{\widetilde h}^2=\int_X\sum_i
|\chi_i\sigma_i|_{h_{H_\ell}}^2|\gamma|_{h^k}^2\,dV\ge e^{-bM}\|u\|_{h^k}^2,
\end{equation}
for some constant $b>0$.

Now we compare \eqref{eq64} with \eqref{eq65}. Choose $A$ sufficiently large such that $cA^2>b$ and then choose sufficiently large $M$ and $k$ such that
\begin{equation}\label{eq66}
\sum_i\|(I-\Pi)\mathcal E_i u\|_{\widetilde h}^2
\leq\frac12e^{-bM}\|u\|_{h^k}^2.
\end{equation}
Therefore we can deduce from \eqref{eq65} and \eqref{eq66} that
\[
\begin{split}
\|\mathcal T u\|_{\widetilde h}^2
&=\sum_i\|\Pi\mathcal E_i u\|_{\widetilde h}^2=\sum_i\|\mathcal E_i u\|_{\widetilde h}^2
-\sum_i\|(I-\Pi)\mathcal E_i u\|_{\widetilde h}^2\\
&\geq\frac12e^{-bM}\|u\|_{h^k}^2,
\end{split}
\]
which proves the first inequality in \eqref{eq58} with $a=\frac12e^{-bM}$.

Set $N=\rank\Pi$ and let $e_1,\ldots,e_N$ be an orthonormal basis $\operatorname{Ran}\Pi$ with respect to $\widetilde h$ For a unit harmonic form $u$, the Cauchy-Schwarz inequality and Corollary \ref{cor2} yield
\begin{align*}
 |\langle\mathcal E_i u,e_j\rangle_{\widetilde h}|^2
 &\le\left(\int_{D_i}|\gamma|_{h^k}^2\,dV\right)
      \left(\int_{D_i}\chi_i^2|\sigma_i|_{h_{H_\ell}}^2|e_j|_{\widetilde h}^2\,dV\right)
 \\
 &\le C_A r^{2q}\int_{D_i}\chi_i^2|\sigma_i|_{h_{H_\ell}}^2|e_j|_{\widetilde h}^2\,dV.
\end{align*}
Let $\mathcal E_i^*$ denote the Hilbert adjoint of $\mathcal E_i$ with respect to $h^k$. By taking the supremum over all such $u$ we get
\begin{equation}\label{eq67}
\|\mathcal E_i^*e_j\|_{h^k}^2
\le C_A r^{2q}\int_{D_i}\chi_i^2|\sigma_i|_{h_{H_\ell}}^2|e_j|_{\widetilde h}^2\,dV,
\end{equation}
We finally deduce from the Bessel's inequality, \eqref{eq67} and Lemma \ref{lem6} that
\begin{equation}\label{eq68}
\begin{split}
\sum_{\alpha=1}^d\|\mathcal T u_\alpha\|_{\widetilde h}^2&=
\sum_{j=1}^N\sum_i\sum_{\alpha=1}^d
|\langle\mathcal E_i u_\alpha,e_j\rangle_{\widetilde h}|^2=\sum_{j=1}^N\sum_i\sum_{\alpha=1}^d
|\langle u_\alpha,\mathcal E_i^*e_j\rangle_{h^k}|^2\\
&\le \sum_{j=1}^N\sum_i\|\mathcal E_i^*e_j\|_{h^k}^2\le C_A r^{2q}\sum_{j=1}^N\int_X
\big(\sum_i\chi_i^2|\sigma_i|_{h_{H_\ell}}^2\big)|e_j|_{\widetilde h}^2\,dV\\
&\le \widetilde C_A r^{2q}\sum_{j=1}^N\int_X
|e_j|_{\widetilde h}^2\,dV=\widetilde C_A r^{2q}N\\
&\le Ck^{n-q}.
\end{split}
\end{equation}
The proof is completed.
\end{proof}
The factor $k^{-q}=r^{2q}$ is precisely the one inherited from Berndtsson's local mass estimate \cite{Berndtsson}. Now we arrive at the

\begin{proof}[Proof of Theorem \ref{thm1}]
Note that $\I(h^k)\subset\OO_X$, then the $q=0$ case follows from \cite[Theorem~5.1]{DemaillyMorse}. For $1\le q\le n$ and $F=E\otimes K_X^{-1}$, Lemma \ref{lem2} identifies the cohomology group with $\HH_k^{n,q}(F)$. For an orthonormal basis $u_1,\ldots,u_{d_k}$ of $\HH_k^{n,q}(F)$, we have
\[
 ad_k\le\sum_{j=1}^{d_k}\|\mathcal T u_j\|_{\widetilde h}^2\le Ck^{n-q}
\]
by Proposition \ref{prop3}, which completes the proof.
\end{proof}

\begin{proof}[Proof of Corollary \ref{cor1}]
The case $q=0$ is trivial as in the proof of Theorem \ref{thm1}. Consider the case when $1\le q\le n$. By nefness of $L$, for every sufficiently large integer $k$, we can choose a smooth Hermitian metric $g_k$ on $L$ such that
\begin{equation}\label{eq69}
k\,\Theta_{g_k}(L)\ge -\omega,
\end{equation}
with $g_k=h_0e^{-U_k}, U_k\le0$. We repeat the arguments of Sections~\ref{sec3} and \ref{sec4}, with the singular metric $h$ replaced by the smooth metric $g_k$. The curvature assumption \eqref{eq69} ensures that all local estimates in Section~\ref{sec3} hold with constants independent of $k$ and of $g_k$. Note that in \eqref{eq38}, $h_{0,k}=h_Gh_0^kh_{H_\ell}$ and $h_{1,k}=h_Gg_k^kh_{H_\ell}$. By repeating the proof of Lemma \ref{lem5} (cf. \eqref{eq41}), we conclude that for sufficiently large $M$
\[
\Theta(H_{k,\ell}\otimes K_X^{-1})\geq c_0\ell\omega.
\]
Then the results in Section \ref{sec4} hold uniformly in $g_k$. From Proposition \ref{prop3}, we can also show that \eqref{eq58} is independent of $g_k$. The proof is done by repeating the proof of Theorem \ref{thm1}.
\end{proof}

\appendix

\section{A smooth graph approximation for weighted maximal complexes}

Let $\Omega$ be a relatively compact coordinate domain and let $V=\bigoplus_{j=0}^N V^j$ be a smooth Hermitian bundle with orthogonal grading. Fix a smooth positive volume form $dV$. Let $d$ be a smooth first-order differential of degree one such that $d^2=0$. For a measurable function $\psi$, finite almost everywhere, with $e^{-\psi}$ locally bounded below by a positive constant, put
\[
 H_\psi=L^2(\Omega,V;e^{-\psi}),\quad
 \|b\|_\psi^2=\int_\Omega |b|_{V}^2e^{-\psi}\,dV.
\]
The local lower bound embeds $H_\psi$ continuously into $L^2_{\mathrm{loc}}$. It's easy to see that $d$ is a closed operator. Whenever $d_\psi$ is densely defined, $d_\psi^*$ denotes its Hilbert adjoint.

\begin{lemma}\label{lem8}
Let $\phi_\nu$ be smooth real functions on $\Omega$, and let $\phi$ be measurable and finite almost everywhere. Suppose that $d_\phi$ is densely defined and that
\begin{equation}\label{eq70}
 e^{-\phi_\nu}\longrightarrow e^{-\phi}\quad\text{almost everywhere},
 \quad c\le e^{-\phi_\nu}\le Ce^{-\phi},
\end{equation}
where $c,C>0$ are independent of $\nu$. For every compactly supported $v\in\Dom d_\phi\cap\Dom d_\phi^*$ there are $v_\nu\in C_c^\infty(\Omega,V)$, $\supp(v_\nu)=K\Subset \Omega$, such that
\begin{equation}\label{eq71}
 \begin{aligned}
 e^{-\phi_\nu/2}v_\nu&\longrightarrow e^{-\phi/2}v,\\
 e^{-\phi_\nu/2}dv_\nu&\longrightarrow e^{-\phi/2}dv,\\
 e^{-\phi_\nu/2}d_{\phi_\nu}^*v_\nu
     &\longrightarrow e^{-\phi/2}d_\phi^*v
 \end{aligned}
 \quad\text{strongly in }L^2(\Omega,V).
\end{equation}
\end{lemma}

\begin{proof}
Put $H=L^2(\Omega,V)$ with the original smooth metrics and volume form. The equation $\|e^{-\psi/2}b\|_H=\|b\|_\psi$
identifies each weighted space isometrically with $H$. On $H$, we define
\[
 D_\nu a=e^{-\phi_\nu/2}d(e^{\phi_\nu/2}a),\quad
 Da=e^{-\phi/2}d(e^{\phi/2}a),
\]
with domains
\[
 \begin{aligned}
 \Dom D_\nu
   &=\{a\in H:e^{\phi_\nu/2}a\in\Dom d_{\phi_\nu}\},\\
 \Dom D
   &=\{a\in H:e^{\phi/2}a\in\Dom d_\phi\}.
 \end{aligned}
\]
These are closed densely defined differentials. The weighted adjoint pairing gives
\[
 D_\nu^*a=e^{-\phi_\nu/2}d_{\phi_\nu}^*(e^{\phi_\nu/2}a),
 \quad
 D^*a=e^{-\phi/2}d_\phi^*(e^{\phi/2}a),
\]
on the corresponding adjoint domains. We first prove the desired result for $D_\nu,D$ on $H$, and then return to the weighted spaces.

Let $b\in\Dom d_\phi$. The upper bound in \eqref{eq70} ensures that $b$ belongs to $\Dom d_{\phi_\nu}$ for every $\nu$.
Moreover,
\[
 |e^{-\phi_\nu/2}-e^{-\phi/2}|^2(|b|_{V}^2+|db|_{V}^2)
 \le (C+1)e^{-\phi}(|b|_{V}^2+|db|_{V}^2).
\]
The dominated convergence theorem therefore yields
\[
 (e^{-\phi_\nu/2}b,e^{-\phi_\nu/2}db)
 \longrightarrow(e^{-\phi/2}b,e^{-\phi/2}db)
 \quad\text{strongly in }H\oplus H.
\]
This means for all $x\in \graph D$, we can find a sequence $x_\nu\in \graph D_\nu$ such that $x_\nu\to x$ strongly.

Suppose $(a_\nu,b_\nu)\in\operatorname{graph}D_\nu$ converges weakly to $(a,b)$ in $H\oplus H$. Let $d^\dagger$ be the formal adjoint with respect to the fixed smooth metrics and volume form. Then we have
\[
 \langle b_\nu,e^{\phi_\nu/2}\zeta\rangle_H
 =\langle a_\nu,e^{\phi_\nu/2}d^\dagger\zeta\rangle_H.
\]
The lower bound in \eqref{eq70} implies $e^{\phi_\nu/2}\le c^{-1/2}$, so the dominated convergence theorem again implies
\[
 e^{\phi_\nu/2}\xi\longrightarrow e^{\phi/2}\xi
 \quad\text{strongly in }H\quad(\xi\in H).
\]
Hence we obtain $d(e^{\phi/2}a)=e^{\phi/2}b$. Since
\[
 \|e^{\phi/2}a\|_\phi=\|a\|_H,\quad
 \|e^{\phi/2}b\|_\phi=\|b\|_H,
\]
the pair $(a,b)$ belongs to $\operatorname{graph}D$.

Let $P_\nu,P$ be the orthogonal projections onto these graphs and fix $z\in H\oplus H$. Any weak subsequential limit $y$ of the bounded sequence $P_\nu z$ belongs to $\operatorname{graph}D$. For $x$ in this graph, choose strong graph approximations $x_\nu$. The identities $\langle z-P_\nu z,x_\nu\rangle_{H\oplus H}=0$ give $z-y\perp\operatorname{graph}D$, so $y=Pz$. Hence $P_\nu z\rightharpoonup Pz$, and
\[
 \|P_\nu z\|_{H\oplus H}^2=\langle P_\nu z,z\rangle_{H\oplus H}
 \longrightarrow\langle Pz,z\rangle_{H\oplus H}=\|Pz\|_{H\oplus H}^2.
\]
It follows that $P_\nu z\to Pz$ strongly for every $z$.

Write the bounded self-adjoint projection $P$ as
\begin{equation}
 P=\begin{pmatrix}X&Z^*\\ Z&I-Y\end{pmatrix},
 \quad X=X^*,\quad Y=Y^*.
\end{equation}
The definition of the Hilbert adjoint shows
\[
 (\operatorname{graph}D)^\perp
 =\{(-D^*b,b):b\in\Dom D^*\}.
\]
Indeed, orthogonality of $(r,s)$ to all $(x,Dx)$ is exactly $\langle Dx,s\rangle_{H}=\langle x,-r\rangle_{H}$ for all $x\in\Dom D$, which means $s\in\Dom D^*$ and $D^*s=-r$.

For $h\in H$ we have
\[
 \begin{aligned}
 P(h,0)&=(Xh,Zh),&
 P(0,h)&=(Z^*h,(I-Y)h),\\
 (I-P)(h,0)&=((I-X)h,-Zh),&
 (I-P)(0,h)&=(-Z^*h,Yh).
 \end{aligned}
\]
Reading the first row in $\operatorname{graph}D$ and the second row in $(\graph D)^\perp$ yields
\begin{equation}\label{eq73}
 DX=Z,\quad DZ^*=I-Y,\quad
 D^*Z=I-X,\quad D^*Y=Z^*.
\end{equation}
In particular,
\[
 XH,Z^*H\subset\Dom D,\quad
 ZH,YH\subset\Dom D^*.
\]

Since $D$ is a differential, $\Ran D\subset\ker D$. Also $\Ran D^*\subset\ker D^*$. If $h\in\ker D$, then $P(h,0)=(h,0)$, so $Xh=h$. If $h\in\ker D^*$, then $(0,h)$ is orthogonal to the graph, so $P(0,h)=0$ and $Yh=h$. Consequently
\[
 X|_{\ker D}=I,\quad Y|_{\ker D^*}=I.
\]
By \eqref{eq73},
\[
 \Ran(I-Y)\subset\ker D,\quad
 \Ran(I-X)\subset\ker D^*.
\]
It follows that $X(I-Y)=I-Y$, hence $XY=X+Y-I$. By taking adjoints we have $YX=X+Y-I$. Put
\begin{equation}
 W:=X+Y-I=XY=YX.
\end{equation}
$W=XY$ shows that $WH\subset\Dom D$; the representation $W=YX$ shows that $WH\subset\Dom D^*$. Likewise, $Z=DX$ and $ZH\subset\Dom D$, while \eqref{eq73} already yields $ZH\subset\Dom D^*$. The argument for $Z^*$ is symmetric. Thus $W,Z,Z^*$ all map $H$ into $\Dom D^*$.

It follows from \eqref{eq73} also that
\[
 \begin{aligned}
 Z(I-Y)&=DX(I-Y)=D(I-Y)=0,\\
 Z^*(I-X)&=D^*Y(I-X)=D^*(I-X)=0.
 \end{aligned}
\]
Therefore $DW=DXY=ZY=Z$ and $D^*W=D^*YX=Z^*X=Z^*$. We obtain
\begin{equation}\label{eq75}
 \begin{array}{lll}
 DW=Z,&DZ^*=I-Y,&DZ=0,\\
 D^*W=Z^*,&D^*Z=I-X,&D^*Z^*=0.
 \end{array}
\end{equation}
All these  compositions are bounded on $H$. Every argument in Steps~3-4 applies to $D_\nu$ as well. Write $X_\nu,Y_\nu,Z_\nu$ for the corresponding blocks of $P_\nu$ and set $W_\nu=X_\nu+Y_\nu-I$. Strong convergence of $P_\nu$ gives strong convergence of each of its four blocks and therefore also $W_\nu\to W$ strongly.

\smallskip\noindent
\emph{Step 5.} Let $a\in\Dom D\cap\Dom D^*$. Since $P(a,Da)=(a,Da)$ and $P(-D^*a,a)=0$, we get
\[
 \begin{aligned}
 Z^*Da&=(I-X)a,& Za&=YDa,\\
 ZD^*a&=(I-Y)a,& Z^*a&=XD^*a.
 \end{aligned}
\]
In particular, there is the exact decomposition
\[
 Wa+Z^*Da+ZD^*a
 =\bigl[(X+Y-I)+(I-X)+(I-Y)\bigr]a=a.
\]
Hence, we define
\begin{equation}
 a_\nu=W_\nu a+Z_\nu^*Da+Z_\nu D^*a.
\end{equation}
Every term belongs to $\Dom D_\nu\cap\Dom D_\nu^*$. The strong convergence of the above bounded operators implies $a_\nu\to a$. Using \eqref{eq75} with subscripts $\nu$, we also have
\[
 \begin{aligned}
 D_\nu a_\nu
   &=Z_\nu a+(I-Y_\nu)Da
     \longrightarrow Za+(I-Y)Da=Da,\\
 D_\nu^*a_\nu
   &=Z_\nu^*a+(I-X_\nu)D^*a
     \longrightarrow Z^*a+(I-X)D^*a=D^*a.
 \end{aligned}
\]
We have proved
\begin{equation}
 a_\nu\longrightarrow a,\quad
 D_\nu a_\nu\longrightarrow Da,\quad
 D_\nu^*a_\nu\longrightarrow D^*a
 \quad\text{strongly in }H.
\end{equation}

Now we apply the construction to $a=e^{-\phi/2}v$. Since graph projections need not preserve supports, choose a real $\eta\in C_c^\infty(\Omega)$ equal to one near $\supp v$. The commutator $B_\eta=[d,\eta]$ is a bounded zeroth-order operator on $H$. Scalar weight multiplication commutes with it, so
\[
 D_\nu(\eta b)=\eta D_\nu b+B_\eta b
 \quad(b\in\Dom D_\nu).
\]
To verify the corresponding adjoint formula, let $e\in\Dom D_\nu^*$ and test against arbitrary $b\in\Dom D_\nu$. Since $\eta b\in\Dom D_\nu$,
\[
 \langle D_\nu b,\eta e\rangle_H=\langle D_\nu(\eta b)-B_\eta b,e\rangle_H=\langle b,\eta D_\nu^*e-B_\eta^*e\rangle_H.
\]
Consequently, $\eta e\in\Dom D_\nu^*$, and the adjoint satisfies
\[
D_\nu^*(\eta e)=\eta D_\nu^*e-B_\eta^*e.
\]
Here $B_\eta^*$ is the pointwise adjoint for the fixed smooth metrics. Both commutators are independent of the scalar weight. The same argument applies to $D,D^*$, including the singular weight.

Because $\eta a=a$, the product rules imply in $H$
\begin{equation}\label{eq78}
 \begin{aligned}
 \eta a_\nu-a&=\eta(a_\nu-a)\rightarrow 0,\\
 D_\nu(\eta a_\nu)-Da
   &=\eta(D_\nu a_\nu-Da)+B_\eta(a_\nu-a)\rightarrow 0,\\
 D_\nu^*(\eta a_\nu)-D^*a
   &=\eta(D_\nu^*a_\nu-D^*a)-B_\eta^*(a_\nu-a)\rightarrow 0.
 \end{aligned}
\end{equation}

Hence
\[
 \widetilde v_\nu=e^{\phi_\nu/2}\eta a_\nu
\]
belongs to $\Dom d_{\phi_\nu}\cap\Dom d_{\phi_\nu}^*$ and is supported in $\supp\eta$. By the definitions of $D_\nu$ and $D_\nu^*$, the three limits in \eqref{eq71} hold with $\widetilde v_\nu$ in place of $v_\nu$.

Since $\phi_\nu$ is smooth, the formal weighted adjoint $d_\nu^\dagger$ has smooth coefficients. Testing against smooth compactly supported sections shows that it agrees distributionally with $d_{\phi_\nu}^*$ on $\widetilde v_\nu$.
Choose a fixed compact neighborhood $K'\Subset\Omega$ of $\supp\eta$. Interior Friedrichs approximation for the first-order system
\[
 b\longmapsto (db,d_\nu^\dagger b)
\]
shows $b_{\nu,m}\in C_c^\infty(\Omega,V)$ supported in $K'$ such that, as $m\to\infty$,
\[
 b_{\nu,m}\to\widetilde v_\nu,\quad
 db_{\nu,m}\to d\widetilde v_\nu,\quad
 d_\nu^\dagger b_{\nu,m}\to d_{\phi_\nu}^*\widetilde v_\nu
\]
in ordinary $L^2$ \cite[Sections~1-3]{Friedrichs}. 

On $K'$ there are constants $0<m_\nu\le M_\nu<\infty$ with $m_\nu\le e^{-\phi_\nu}\le M_\nu$. Thus ordinary and weighted $L^2$ convergence are equivalent there. Moreover, integration by parts against arbitrary vectors in $\Dom d_{\phi_\nu}$ shows that every $b_{\nu,m}$ belongs to $\Dom d_{\phi_\nu}^*$ and that its Hilbert adjoint equals $d_\nu^\dagger b_{\nu,m}$. Choose $m=m(\nu)$ large enough and set $v_\nu=b_{\nu,m(\nu)}$ so that
\[
 \|v_\nu-\widetilde v_\nu\|_{\phi_\nu}
  +\|dv_\nu-d\widetilde v_\nu\|_{\phi_\nu}
  +\|d_{\phi_\nu}^*v_\nu-d_{\phi_\nu}^*\widetilde v_\nu\|_{\phi_\nu} <\nu^{-1}.
\]
Combining this estimate with \eqref{eq78} proves \eqref{eq71}. The proof is completed.
\end{proof}

\section{The choice of local charts}
Let $B_\omega(x,r)=\{y\in X:d_\omega(x,y)<r\}$ denote the geodesic ball. We use holomorphic coordinate balls with bounded overlap.

\begin{lemma}\label{lem9}
Let $(X,\omega)$ be a compact Hermitian manifold of complex dimension $n$. One can fix, for each $x\in X$, a holomorphic coordinate chart $z_x$ centered at $x$ and a number $r_0>0$ such that
\[
 z_x(x)=0,\quad
 B(x,r):=z_x^{-1}\bigl(\{z\in\C^n:|z|<r\}\bigr)
 \quad(0<r<r_0),
\]
and the following statements hold.
\begin{enumerate}[label=\textup{(\roman*)}]
\item The charts are defined on a uniform coordinate radius, satisfy $\omega(x)=\beta$. After shrinking $r_0$ if necessary,
\[
 B_\omega(x,r/2)\subset B(x,r)\subset B_\omega(x,2r),
 \quad 0<r<r_0,
\]
\item For every sufficiently small $r>0$, there are centers $x_1,\ldots,x_{N(r)}$ with
\[
 d_\omega(x_i,x_j)\ge r/4\quad(i\ne j),\quad
 X=\bigcup_{i=1}^{N(r)}B(x_i,r).
\]
\item For each fixed $A\ge1$, there are $r_A>0$ and $N_A<\infty$ independent of $r$ such that, for $0<r<r_A$,
\[
 \sum_{i=1}^{N(r)}\one_{B(x_i,Ar)}(y)\le N_A
 \quad(y\in X).
\]
\end{enumerate}

\end{lemma}

\begin{proof}
For (i), choose finitely many holomorphic coordinate charts $z_j:U_j\longrightarrow\Omega_j\subset\C^n$ and open sets $V_j\Subset U_j$ such that $X=\bigcup_{j=1}^J V_j$. For every $x\in X$, define
\[
 j(x)=\min\{j:x\in V_j\}.
\]
Let $j=j(x)$ and write the metric at $x$ in the coordinates $z_j$ as a positive Hermitian form with matrix $G_j(x)$. Choose a complex linear map $A_x$ that sends this Hermitian form isometrically to the standard Hermitian form, equivalently satisfying $A_x^*A_x=G_j(x)$ in the corresponding vector convention. Set
\[
 z_x(y)=A_x\bigl(z_j(y)-z_j(x)\bigr).
\]
Thus $z_x(x)=0$ and the coefficient matrix of $\omega$ at the origin is the identity. The finite collection of coordinate charts ensures that the estimates are uniform in $x$, even when the chosen chart varies.

Since the atlas is finite and $\overline{V_j}\subset U_j$. With smooth positive definiteness of $\omega$ on the finitely many $\overline{V_j}$ we obtain uniform bounds for $\|A_x\|$ and $\|A_x^{-1}\|$. Hence the $z_x$ are all defined on a common Euclidean ball. In particular, for a uniformly small $r_0$,
\[
 \tfrac12\beta\le\omega\le
 2\beta\quad\text{on }B(x,r_0)
\]
for every $x\in X$. Then $B(x,r)\subset B_\omega(x,2r)$. Conversely, any path starting at $x$ and exiting $B(x,r)$ has length at least $r/\sqrt{2}$.
Thus $B_\omega(x,r/2)\subset B(x,r)$, proving (i).

For the second one, fix $r>0$ small and choose a maximal finite $(r/4)$-separated subset $\{x_1,\ldots,x_{N(r)}\}$ in the compact metric space $(X,d_\omega)$. Here an $r/4$-separated set is a subset of $X$ whose distinct points are at least $r/4$ apart with respect to the metric \(d_\omega\).
Maximality implies $X=\bigcup_i B_\omega(x_i,r/4)$: otherwise a point outside these balls could be added to the separated set. Then by (i) we have $X=\bigcup_i B(x_i,r)$, proving (ii).

To prove (iii), note that the geodesic balls $B_\omega(x_i,r/12)$ are pairwise disjoint. By compactness and the small-ball volume comparison, there are $v_-,v_+>0$ and $r_{\mathrm{vol}}>0$ such that
\[
 v_-r^{2n}\le\operatorname{Vol}_\omega B_\omega(x,r)
 \le v_+r^{2n}\quad(x\in X,\ 0<r<r_{\mathrm{vol}}).
\]
Fix $A\ge1$ and $y\in X$. If $y\in B(x_i,Ar)$, then by (i) $d_\omega(y,x_i)<2Ar$. The pairwise disjoint balls $B_\omega(x_i,r/12)$ associated with these indices all lie in $B_\omega(y,(2A+1/12)r)$. For $r$ small enough that $(2A+1/12)r<r_{\mathrm{vol}}$ and $Ar<r_0$, the volume comparison yields
\[
 \#\{i:y\in B(x_i,Ar)\}\,v_-(r/12)^{2n}
 \le v_+\bigl((2A+1/12)r\bigr)^{2n}.
\]
Hence the left-hand cardinality is at most $(v_+/v_-)(24A+1)^{2n}$, proving (iii).

\end{proof}

\section{An identity for non-smooth radial functions}\label{sec8}
Work in the coordinate chart of Proposition~\ref{prop2}, with $\varrho=|z|^2$. In \eqref{eq24}, $\langle S,dd^cf_r\rangle=\langle dd^cS,f_r\rangle$ holds even when $f_r$ is not smooth. We clarify this for the reader's convenience. Recall that $dd^cS+C_0\mu\ge0$ in Proposition \ref{prop2}, which is a positive Radon measure. Then $dd^cS$ is a locally finite signed measure and $\langle dd^cS,f_r\rangle$ makes sense. Let $g_r(s)=(r^2-s)_+$ similarly defined as $f_r(z)$. Choose a nonnegative smooth mollifier $\eta_\varepsilon$ on $\R$, supported in $(-\varepsilon,\varepsilon)$, with $\int_{\R}\eta_\varepsilon(s)ds=1$. Set
\[
g_{r,\varepsilon}:=g_r*\eta_\varepsilon,\quad
f_{r,\varepsilon}(z):=g_{r,\varepsilon}(\varrho(z))\to f_r(z) \ \text {uniformly}.
\]
Then
\begin{equation}\label{eq79}
\langle dd^cS,f_{r,\varepsilon}\rangle
=2\int g_{r,\varepsilon}'(\varrho)\,S\wedge\beta+
\int \eta_\varepsilon(\varrho-r^2)
S\wedge d\varrho\wedge d^c\varrho.
\end{equation}

Write
\[
S\wedge d\varrho\wedge d^c\varrho=b(z)\,d\lambda(z),
\]
where $b\ge0$ is locally integrable. By the coarea formula \cite{Hor, PSWW}, there is a locally integrable function $A(s)$ satisfying
\[
\int \eta_\varepsilon(\varrho-r^2)b(z)\,d\lambda(z)=
\int \eta_\varepsilon(s-r^2)A(s)\,ds.
\]
Then
\begin{equation}\label{eq80}
\int\delta(\varrho-r^2)S\wedge d\varrho\wedge d^c\varrho=A(r^2)
\end{equation}
holds for a.e. $r>0$. Recall that $dd^cS$ is a signed measure, by the dominated convergence theorem and \eqref{eq80}, we derive the equalities in \eqref{eq24} by letting $\varepsilon\to 0$ in \eqref{eq79}.

\phantomsection
\addcontentsline{toc}{section}{References}
\enlargethispage{4\baselineskip}

\end{document}